\documentclass[12pt]{article}
\pdfoutput=0
\usepackage{hyperref}
\usepackage{epsfig}
\usepackage{verbatim}
\newcommand{\nocopyright}{
No Copyright\thanks{
The authors hereby waive all copyright
and related or neighboring rights to this work,
and dedicate it to the public domain.
This applies worldwide.
}}
\title{Geometry and the Imagination in Minneapolis}
\author{John H. Conway \and Peter G. Doyle \and Jane Gilman \and William P. Thurston}
\date{
June 1991\\
Version 2.0, 8 April 2018\\
\nocopyright
}
\newcommand{\fig}[3]{
\begin{figure}
\centerline{\mbox{\includegraphics[width=370pt]{figures/#1.eps}}}
\caption{#3}
\label{#2}
\end{figure}
}
\newcommand{\figsize}[4]{
\begin{figure}
\centerline{\mbox{\includegraphics[width=#1]{figures/#2.eps}}}
\caption{#4}
\label{#3}
\end{figure}
}

\newcommand{\dt}[1]{{\em #1}}
\begin{document}

\maketitle

\begin{abstract}
This document consists of the collection of handouts for a
two-week summer workshop entitled 'Geometry and the Imagination',
led by John Conway, Peter Doyle, Jane Gilman and Bill Thurston at
the Geometry Center in Minneapolis, June 17-28, 1991.
The workshop was based on a course `Geometry and the Imagination'
which we had taught twice before at Princeton.
\end{abstract}

\clearpage

\newcommand{\poster}[1]{\includegraphics[width=400pt]{figures/#1.ps}}
\poster{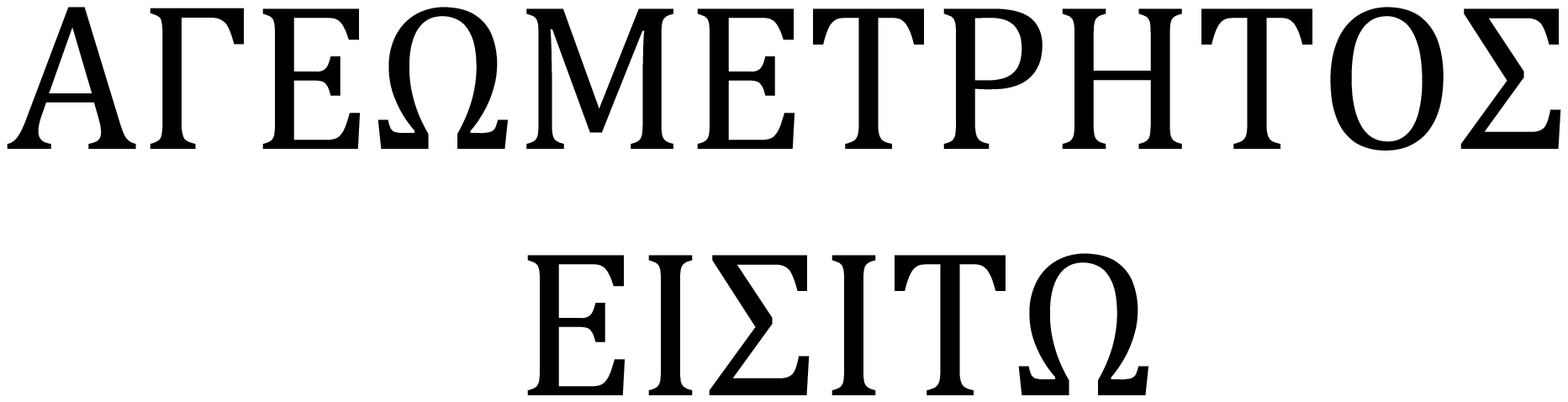}
\poster{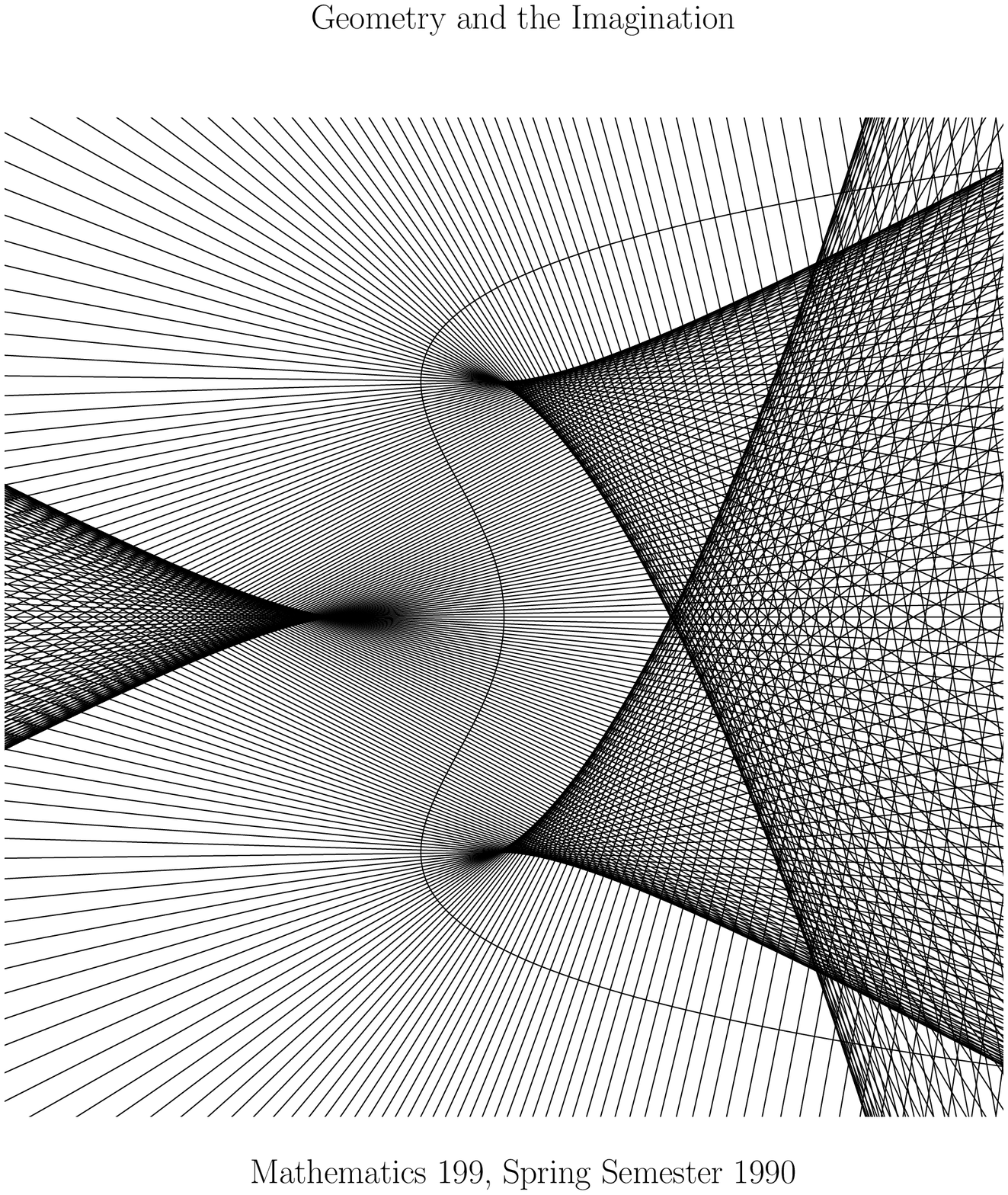}
\poster{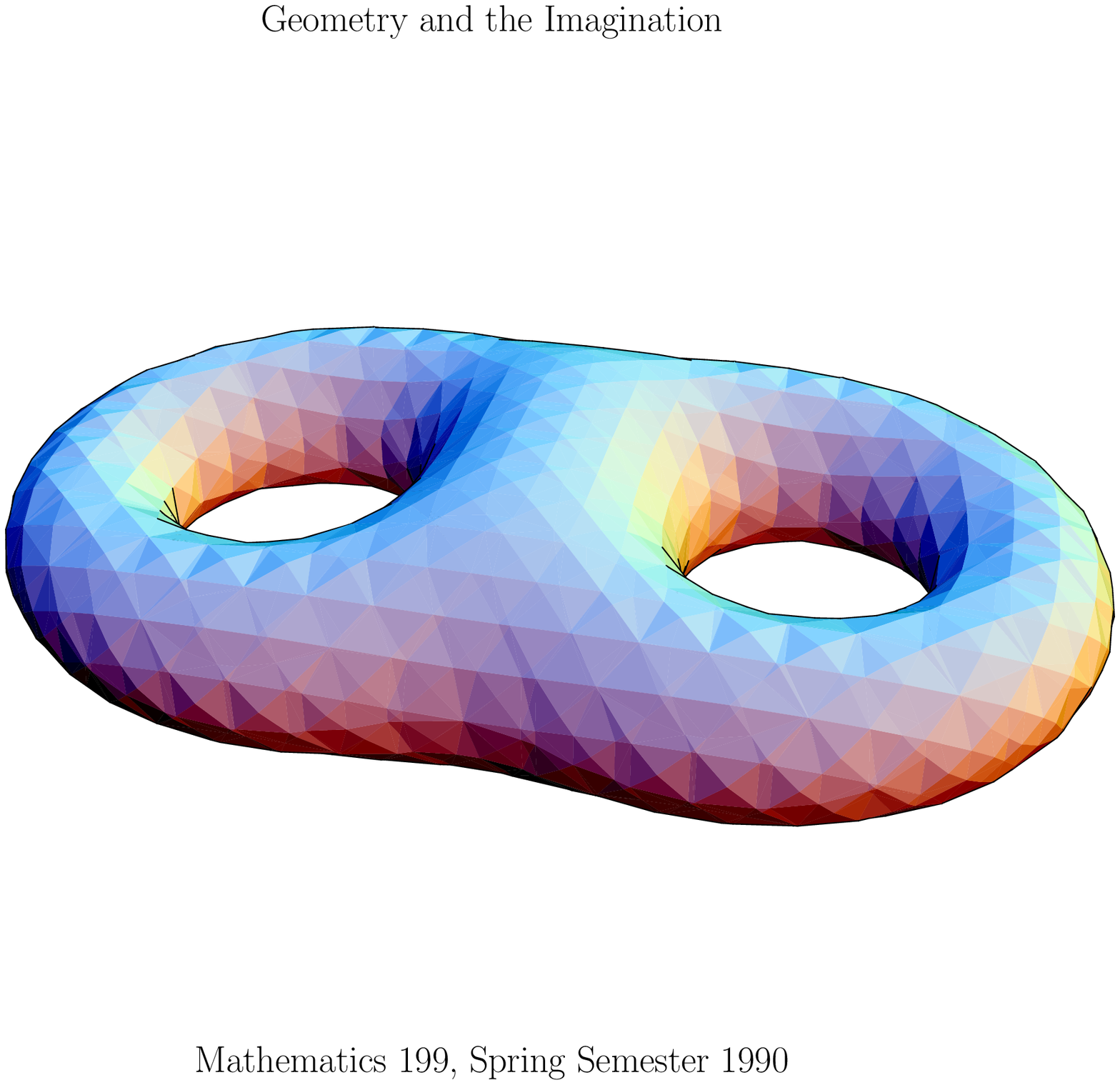}
\poster{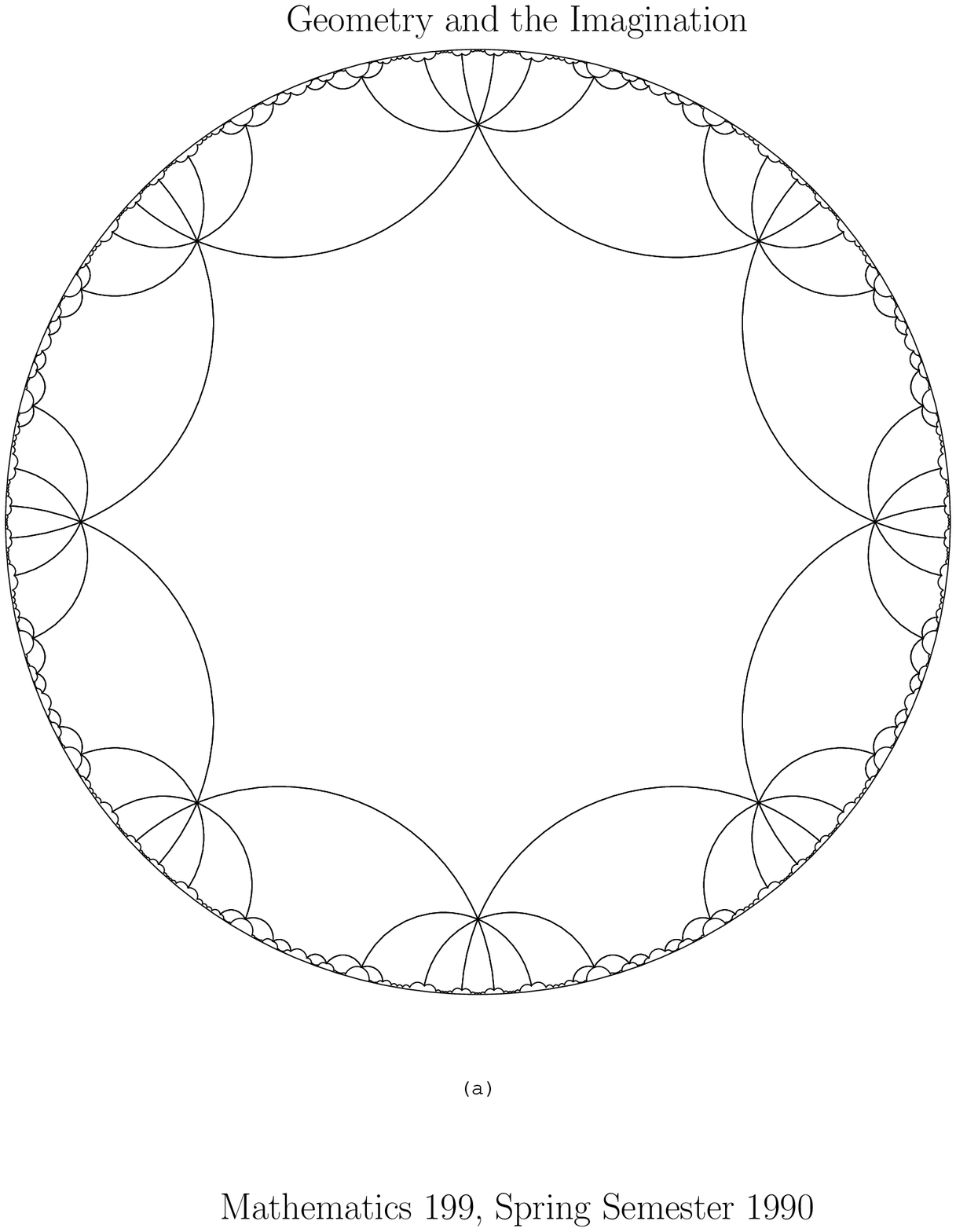}
\poster{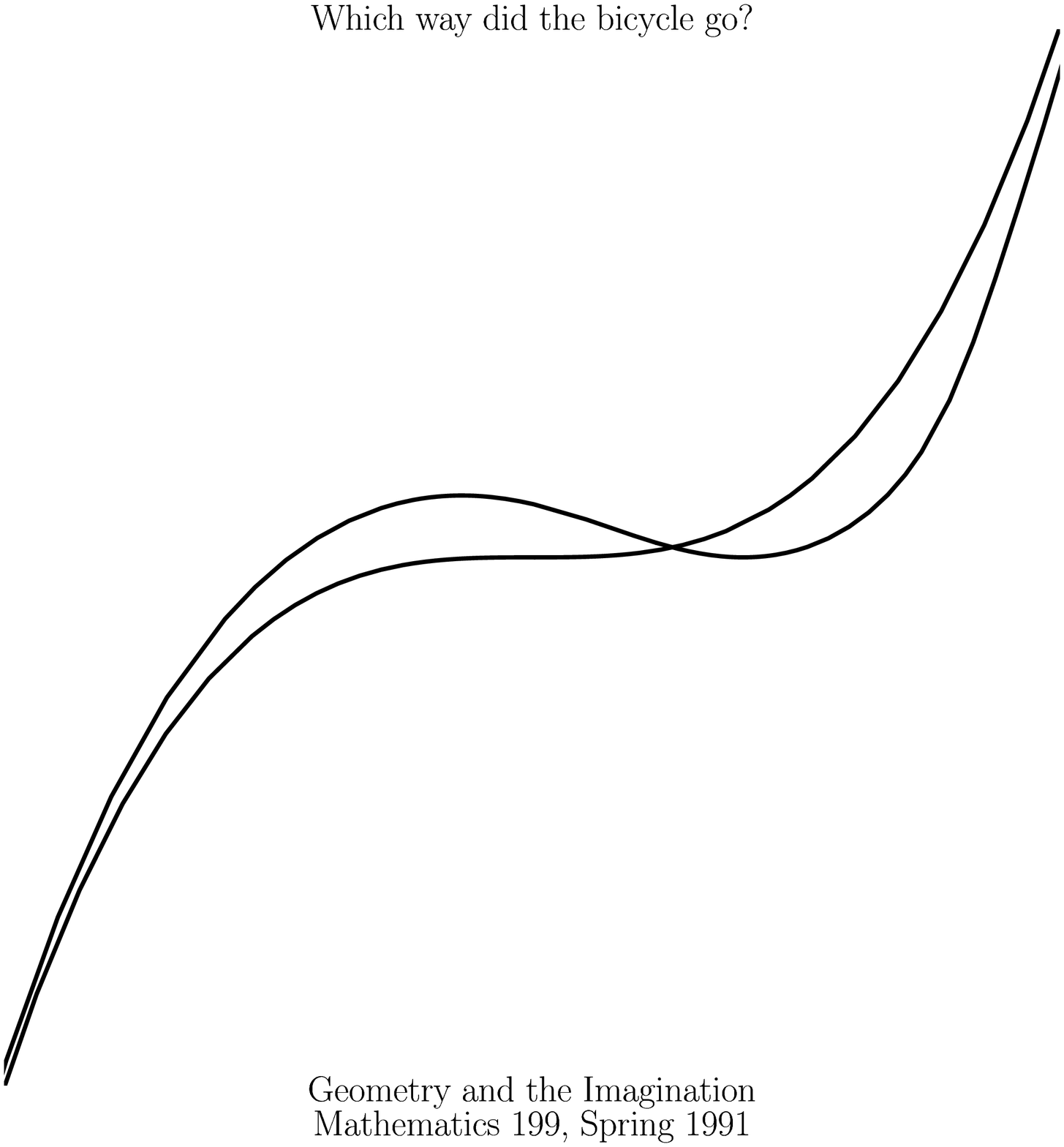}
\poster{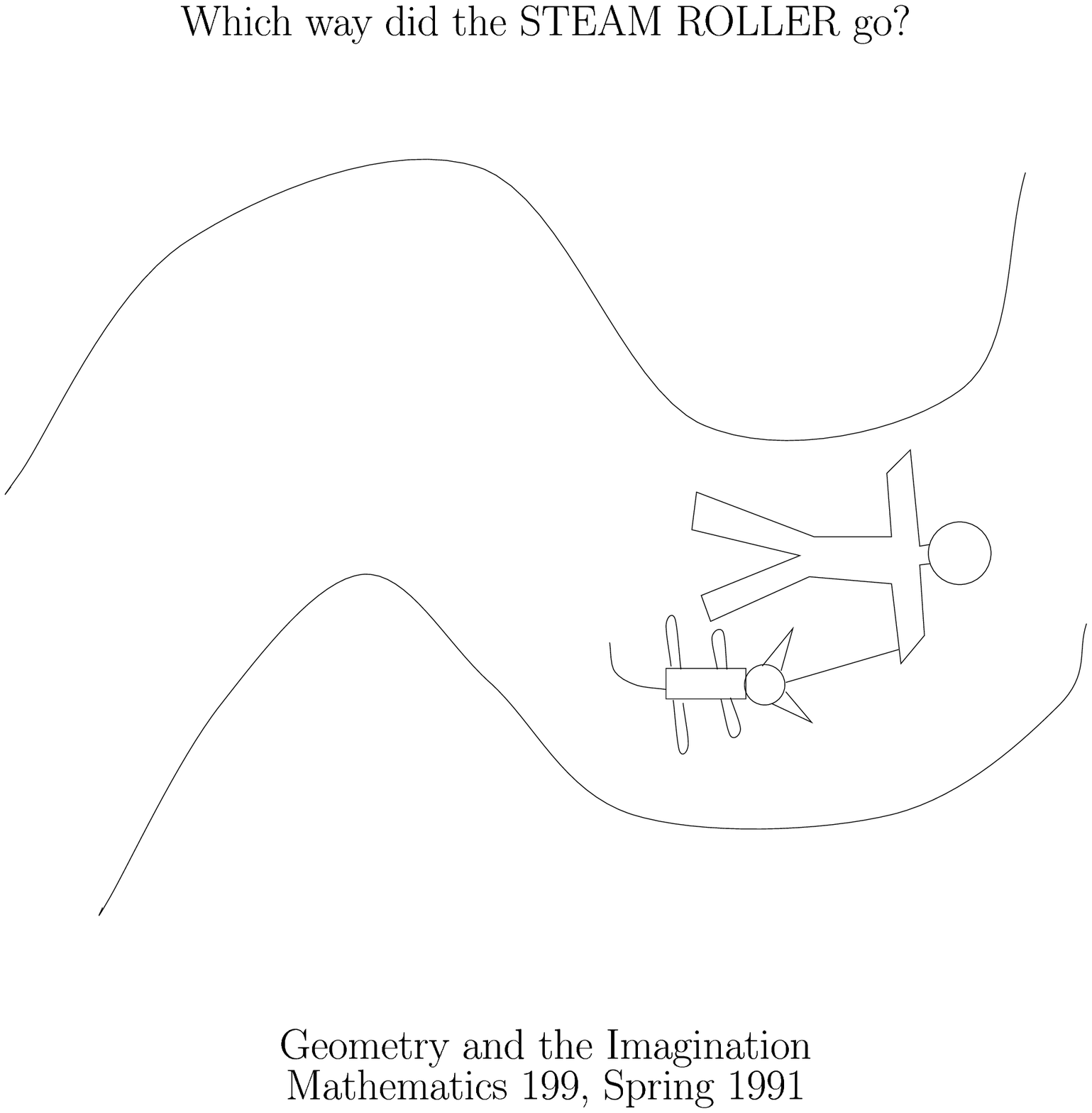}

\clearpage

\section{Preface}

This document consists of the collection of handouts for a
two-week summer workshop entitled 'Geometry and the Imagination',
led by John Conway, Peter Doyle, Jane Gilman and Bill Thurston at
the Geometry Center in Minneapolis, June 17-28, 1991.
The workshop was based on a course `Geometry and the Imagination'
which we had taught twice before at Princeton.

The handouts
do not give a uniform treatment of the topics covered in the workshop: some
ideas were treated almost entirely in class by lecture and
discussion, and other ideas which are fairly extensively documented
were only lightly treated in class.  The motivation for the handouts
was mainly to supplement the class, not to document it.

The primary outside reading was `The Shape of Space', by Jeff Weeks.
Some of the topics discussed in the course which are omitted or only
lightly covered in the handouts are developed well in that book:
in particular, the concepts of extrinsic versus intrinsic topology and
geometry, and two and three dimensional manifolds.  Our approach
to curvature is only partly documented in the handouts.  Activities
with scissors, cabbage, kale, flashlights, polydrons, sewing,
and polyhedra were really live rather than written.

The mix of students---high school students, college undergraduates,
high school teachers and college teachers---was unusual,
and the mode of running a class with the four of us teaching was also unusual.
The mixture of people helped create the tremendous flow of energy
and enthusiasm during the workshop.

Besides the four teachers and the official students, there were
many people who put a lot in to help organize or operate
the course, including Jennifer Alsted, Phil Carlson, Anthony Iano-Fletcher,
Maria Iano-Fletcher, Kathy Gilder, Harvey Keynes, Al Marden, Delle Maxwell,
Jeff Ondich, Tony Phillips, John Sullivan, Margaret Thurston, Angie Vail,
Stan Wagon.

\clearpage

\section{Philosophy}

Welcome to Geometry and the Imagination!

This course aims to convey the richness, diversity, connectedness,
depth and pleasure of mathematics. The title is taken from the classic
book by Hilbert and Cohn-Vossen, ``Geometry and the Imagination'. {\it
Geometry} is taken in a broad sense, as used by mathematicians, to
include such fields as topology and differential geometry as well as
more classical geometry. {\it Imagination}, an essential part of
mathematics, means not only the facility which is imaginative, but also
the facility which calls to mind and manipulates mental images. One aim
of the course is to develop the imagination.

While the mathematical content of the course will be high, we will try
to make it as independent of prior background as possible. Calculus,
for example, is not a prerequisite.

We will emphasize the {\it process} of thinking about mathematics.
Assignments will involve thinking and writing, not just grinding
through formulas.  There will be a strong emphasis on projects and
discussions rather than lectures. All students are expected to get
involved in discussions, within class and without. A  Geometry Room 
on the fifth floor will be reserved for students in the course. The
room will accrete mathematical models, materials for building models,
references related to geometry, questions, responses and (most
important) people.  There will be computer workstations in or near the
geometry room.  You are encouraged to spend your afternoons on the
fifth floor.

The spirit of mathematics is not captured by spending 3 hours solving
20 look-alike homework problems.  Mathematics is thinking, comparing,
analyzing, inventing, and understanding.  The main point is not
quantity or speed---the main point is quality of thought.  The goal is
to reach a
more complete and a better understanding.  We will use materials such
as mirrors, Polydrons, scissors and tissue paper not because we
think this is easier than  solving algebraic equations and differential
equations, but because we think that this is the way to bring thinking
and reasoning to the course.

We are very enthusiastic about this course, and we have many plans to
facilitate your taking charge and learning. While you won't need a
heavy formal background for the course,  you do need  a commitment of
time and energy.

\section{ Organization }

\subsection{People}
We are experimenting with a diverse group of participants in this course:
high school students, high school teachers, college students,
college teachers, and others.

Topics in mathematics often have many levels of meaning, and we hope and
expect that despite---no, because of---the diversity, there will be a lot
for everyone (including
we the staff) to get from the course.  As you think about something,
you come to understand it from different angles, and on successively deeper
levels.

We want to encourage interactions
between all the participants in the course.  It can be quite interesting
for people with sophisticated backgrounds and with elementary backgrounds
to discuss a topic with each other, and the communication can have a high value
in both directions. 

\subsection{Scheduled meetings} 
The officially scheduled morning sessions, which run from 9:00 to 12:30
with a half-hour break in the middle, form the core of the course.
During these sessions, various
kinds of activities will take place. There will be some more-or-less
traditional presentations, but the main emphasis will be on encouraging
you to discover things for yourself. Thus the class will frequently
break into small groups of about 5--7 people for discussions of various
topics.

\subsection{Discussion groups} We want to enable everyone to be engaged
in discussions while at the same time preserving the unity of the
course.  From time to time, we will break into discussion groups of
5--7 people.

Every member of each group is expected to take part in the discussion
and  to make sure that {\it everyone} is involved:  that everyone is
being heard, everyone is listening, that the discussion is not
dominated by one or two people, that everyone understands what is going
on, and that the group sticks to the subject and really digs in.

Each group will have a reporter. The reporters will rotate so that
everyone will serve as reporter during the next two weeks.  The main
role of the reporters during group discussions is to listen, rather than speak.
The reporters should make sure they understand and write down
the key points and ideas from the discussion, and be prepared to summarize
and explain them to the whole class.  

After a suitable time, we will ask for reports to the entire class.
These will not be formal reports.  Rather, we will hold a summary
discussion among the reporters and teachers, with occasional
contributions from others.

\subsection{Texts} The required texts for the course are:
 Weeks, {\it The Shape of Space} and Coxeter, {\it Introduction to
Geometry}.  There are  available at the University Bookstore.

Coxeter's book will mainly be used as a reference book for the course,
but it is also a book that should be useful to you in the future.

Here is a list of reading assignments from
{\it The Shape of Space} by Weeks.
As Weeks suggests it is important to {\it  ``\ldots read slowly and give things
plenty of time to digest''},
as much as is possible in a condensed course of this type.

\begin{itemize}
\item Monday, June 17: Chapters 1 and 2.
\item Tuesday, June 18: Chapter 3.
\item Wednesday, June 19: Chapter 4.
\item Thursday, June 20: Chapter 5, pages 67-77, and Chapter 6, 85-90.

\item Friday, June 21 -- Sunday, June 23: Chapters 7 \& 8.

\item Monday, June 24: Chapters 9 and 10.
\item Tuesday, June 25: Chapter 11 and 12.
\item Wednesday, June 26: Chapter 13.
\item Thursday, June 27: Chapter 16.
\end{itemize}

In addition, there is a long list of  recommended reading. The geometry
room has a small collection of additional books, which you may read
there. There are several copies of some books  which we highly
recommend such as {\it Flatland} by Abbott and {\it What is
Mathematics} by Courant and Robbins. There are single copies of other
books.

\subsection{Other materials} We will be doing a lot of constructions
during class.  Beginning this Tuesday (June 18th), you should bring
with you to class each time:  scissors, tape, ruler, compass, sharp
pencils, plain white paper.  It would be a capital idea to bring extras
to rent to your classmates.

\subsection{Journals} Each participant should keep a journal
for the course.  While assignments given at class meetings go in the journal,
the journal is for much more: for independent questions, ideas,
and projects. The journal is not for class notes,
but for work outside of class.  The style of the journal will vary from
person to person.  Some will find it useful to write short summaries of
what went on in class. Any questions suggested by the class work should
be in the journal.  The questions can be either speculative questions
or more technical questions.  You may also want to write about the nature
of the class meetings and group discussions: what works for you and what
doesn't work, {\it etc.} 

You are encouraged to cooperate with each other in working on anything
in the course, but what you put in your journal should be you.  If it
is something that has emerged from work with other people, write down
who you have worked with.  Ideas that come from other people should be
given proper attribution.  If you have referred to sources other than
the texts for the course, cite them.

Exposition is important.
If you are presenting the solution to a problem, explain what the problem is.
If you are giving an argument, explain what the point is before you
launch into it.  What you should aim for is something that could communicate
to a friend or a colleague a coherent idea of what
you have been thinking and doing in the course.

Your journal should be kept on loose leaf paper. Journals will be
collected every few days and read and commented upon
by the instructors.
If they are on loose
leaf paper, you can hand in those parts which have not yet been read,
and continue to work on further entries. Pages should be numbered
consecutively and except when otherwise instructed, you should hand in
only those pages which have not previously been read. Write your name
on each page, and, in the upper right hand corner of the first page you
hand in each time, list the pages you have handed in (e.g. [7,12] on
page 7 will indicate that you have handed in 6 pages numbered seven to
twelve).

\medskip

Mainly, the journal is for {\it you}.  In addition,
the journals are an important tool by
which we keep in touch with you and what you are thinking about.
Our experience is that it is really fun and enjoyable when someone
lets us into their head.
No matter what your status in this course, 
keep a journal.

Journals will be collected and read as follows:
\begin{itemize}
\item Wed. June 19th
\item Friday June 21st
\item Tuesday, June 25th
\item Thursday June 26th
\end{itemize}
Your entire journal
should be handed in on Friday June 27th with your final project.
We will return final journals by mail.

\subsection{Constructions} Geometry lends itself to constructions and
models, and we will expect everyone to be engaged in model-making.
There will be minor constructions that may take only half an hour and
that everyone does, but we will also expect larger constructions that
may take longer.

\subsection{Final project} We will not have a final exam for the
course, but in its place, you will undertake a major project.  The
major project may be a paper investigating more deeply some topic we
touch on lightly in class.  Alternatively, it might be based on a major
model project, or it might be a computer-based project. To give you
some ideas, a  list of possible projects will be circulated. However,
you are also encouraged to come up with your own ideas for projects. If
possible, your project should have some visual component, for we will
display all of the projects at the end of the course at
the {\it Geometry Fair}. The project will be due on the morning of
Friday June 28th. The fair will be in the afternoon.

\subsection{Geometry room/area}
The fifth floor houses the Geometry Room. We hope that it will actually
spill out into the hallways and corridors and thus become the geometry
{\it area}. Thus the fifth floor will serve as a work and play room for
this course.  This is where you can find mathematical toys, games, models,
displays and construction materials.  Copies of handouts and 
books and other written materials of interest to students in the
course will be kept here as well.  It should also
serve as a place to go if you want to talk to other students in the
course, or to one of the teachers.  Our current plan is to have this
area open from 1:30 to 4:00 PM Monday through Friday,
beginning right away. There
will be a {\it tour} of the area at the end of Monday's morning session.

\section{Bicycle tracks}

Here is a passage from a Sherlock Holmes story, {\it The Adventure of the Priory School} (by Arthur Conan Doyle):

`This track, as you perceive,
was made by a rider who was going from the direction of the school.'

`Or towards it?'

`No, no, my dear Watson.
The more deeply sunk impression is, of course,
the hind wheel, upon which the weight rests.
You perceive several places where it has passed across and obliterated
the more shallow mark of the front one.
It was undoubtedly heading away from the school.'

\begin{enumerate}
\item Discuss the passage above.  
\item Visualize, discuss, and sketch what bicycle tracks look like.
\item When we present actual bicycle tracks, determine
the direction of motion.
\item What else can you tell about the bike from the tracks?
\end{enumerate}

\section{Polyhedra}

A \dt{polyhedron} is the three-dimensional version of a polygon:
it is a chunk of space with flat walls.  In other words, it is
a three-dimensional figure made by gluing polygons together. 
The word is Greek in origin, meaning many-seated.  The plural is
polyhedra.  The polygonal sides of a polyhedron are called its \dt{faces}.

\subsection{Discussion}
Collect some triangles, either the snap-together plastic polydrons or paper
triangles.  Try gluing them together in various ways to form polyhedra.
\begin{enumerate}
\item
Fasten three triangles together at a vertex.  Complete the figure by
adding one more triangle.  Notice how there are three triangles at
{\it every} vertex.  This figure is called a \dt{tetrahedron} because
it has four faces (see the table of Greek number prefixes.)
\item
Fasten triangles together so there are four at every vertex.
How many faces does it have?  From the table of prefixes below, 
deduce its name.

\item
Do the same, with five at each vertex.

\item
What happens when you fasten triangles six per vertex?
\item
What happens when you fasten triangles seven per vertex?
\end{enumerate}

\begin{table}
\begin{center}
\begin{tabular}{r|l}
1&mono\\
2&di\\
3&tri\\
4&tetra\\
5&penta\\
6&hexa\\
7&hepta\\
8&octa\\
9&ennia\\
10&deca\\
11&hendeca\\
12&dodeca\\
13&triskaideca\\
14&tetrakaideca\\
15&pentakaideca \\
16&hexakaideca\\
17&heptakaideca\\
18&octakaideca\\
19&enniakaideca\\
20&icosa \\
\end{tabular}
\end{center}
\caption{The first 20 Greek number prefixes}
\end{table}

\subsection{Homework}
A \dt{regular polygon} is a polygon with all its edges equal and
all angles equal.  A \dt{regular polyhedron} is one whose faces are
regular polygons, all congruent, and having the same number of polygons
at each vertex. 

For homework, construct models of all possible regular polyhedra,
by trying what happens when you fasten together
regular polygons with 3, 4, 5, 6, 7,
{\it etc} sides so the same number come together at each vertex.

Make a table listing the number of faces, vertices, and edges of each.

What should they be called?

\section{Knots}

A mathematical knot is a knotted loop.   For example,
you might take an extension
cord from a drawer and plug one end into the other: this makes
a mathematical knot.

 It might or might not be possible to unknot
it without unplugging the cord.  A knot which can be unknotted is
called an \dt{unknot}.  

Two knots are
considered equivalent if it is possible to rearrange one to
the form of the other, without cutting the loop and without
allowing it to pass through itself.
The reason for using loops of string in the mathematical
definition is that knots in a length of string can always
be undone by pulling the ends through, so any two lengths of string are
equivalent in this sense.

If you drop a knotted loop of string on a table, it crosses over
itself in a certain number of places.  Possibly, there are ways to rearrange
it with fewer crossings---the minimum possible number of crossings
is the \dt{crossing number} of the knot.

\fig{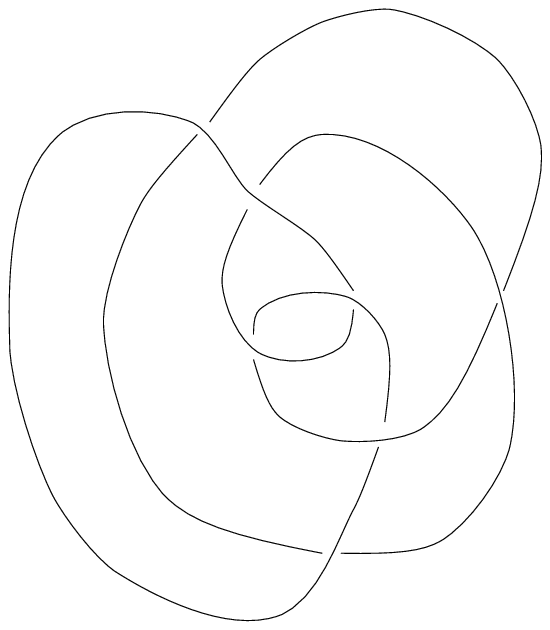}{Drawing of a knot}{This is drawing of a knot has 7
crossings.  Is it possible to rearrange it to have fewer crossings?}{}

\subsection{Discussion}
Make drawings and use short lengths of string to investigate simple knots:
\begin{enumerate}
\item
Are there any knots with one or two crossings?  Why?
\item
How many inequivalent knots are there with three crossings?
\item
How many knots are there with four crossings?
\item
How many knots can you find with five crossings?
\item
How many knots can you find with six crossings?
\end{enumerate}

\section{Maps}

A \dt{map} in the plane
is a collection of vertices and edges (possibly curved) joining the vertices
such that if you cut along the edges the plane falls apart into polygons.
These polygons are called the faces.
A map on the sphere or any other surface is defined similarly.
Two maps are considered to be the same if you can get from one to the
other by a continouous motion of the whole plane.
Thus the two maps in figure \ref{two equivalent maps} are considered to be the same.

\fig{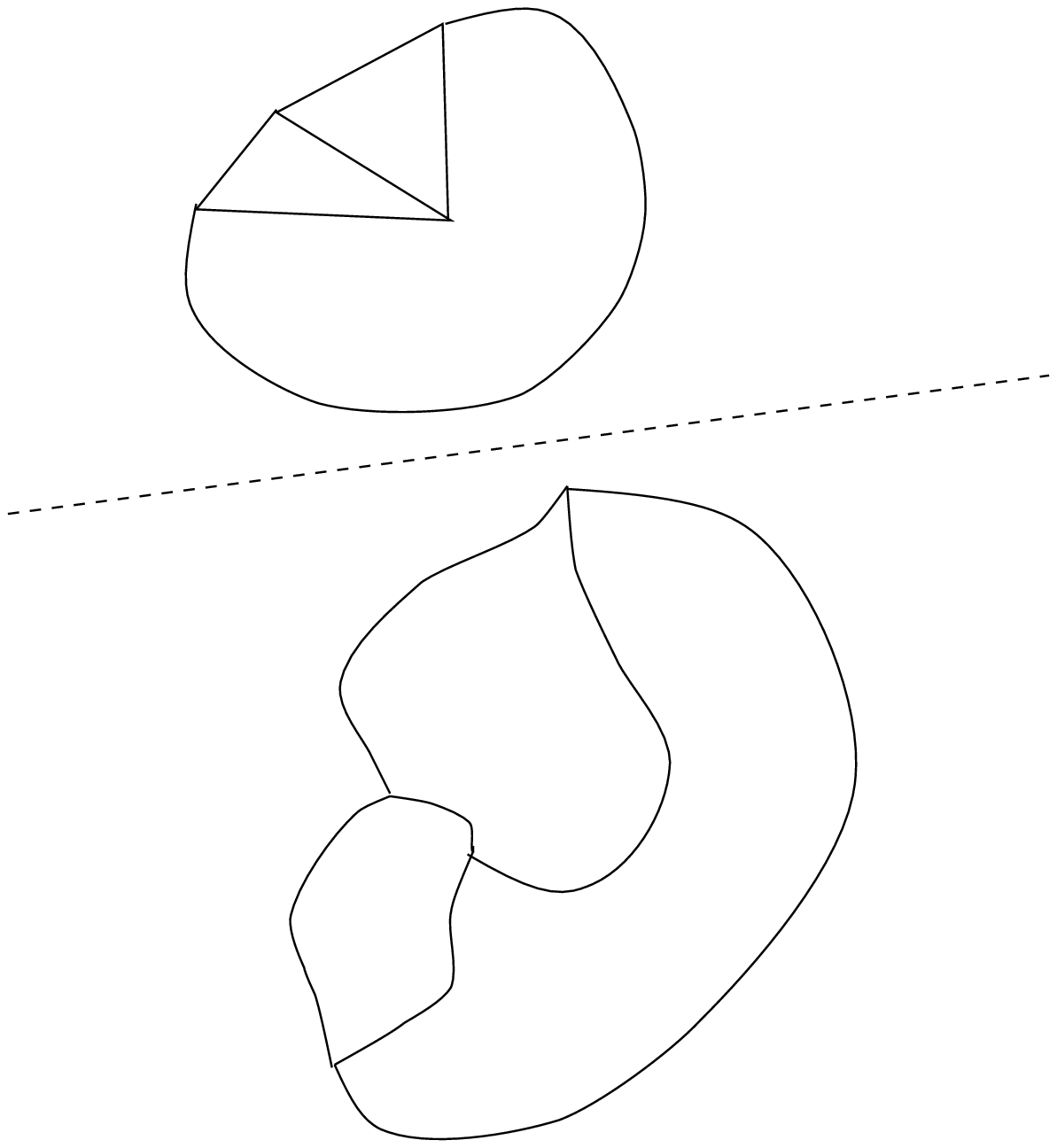}{two equivalent maps}{These two maps are considered
the same (topologically equivalent), because it is possible
to continuously move one to obtain the other.}

A map on the sphere can be represented by a map in the plane by removing
a point from the sphere and then stretching the rest of the sphere out
to cover the plane.
(Imagine popping a balloon and stretching the rubber out onto on the plane,
making sure to stretch the material near the puncture all the way out to
infinity.)

\fig{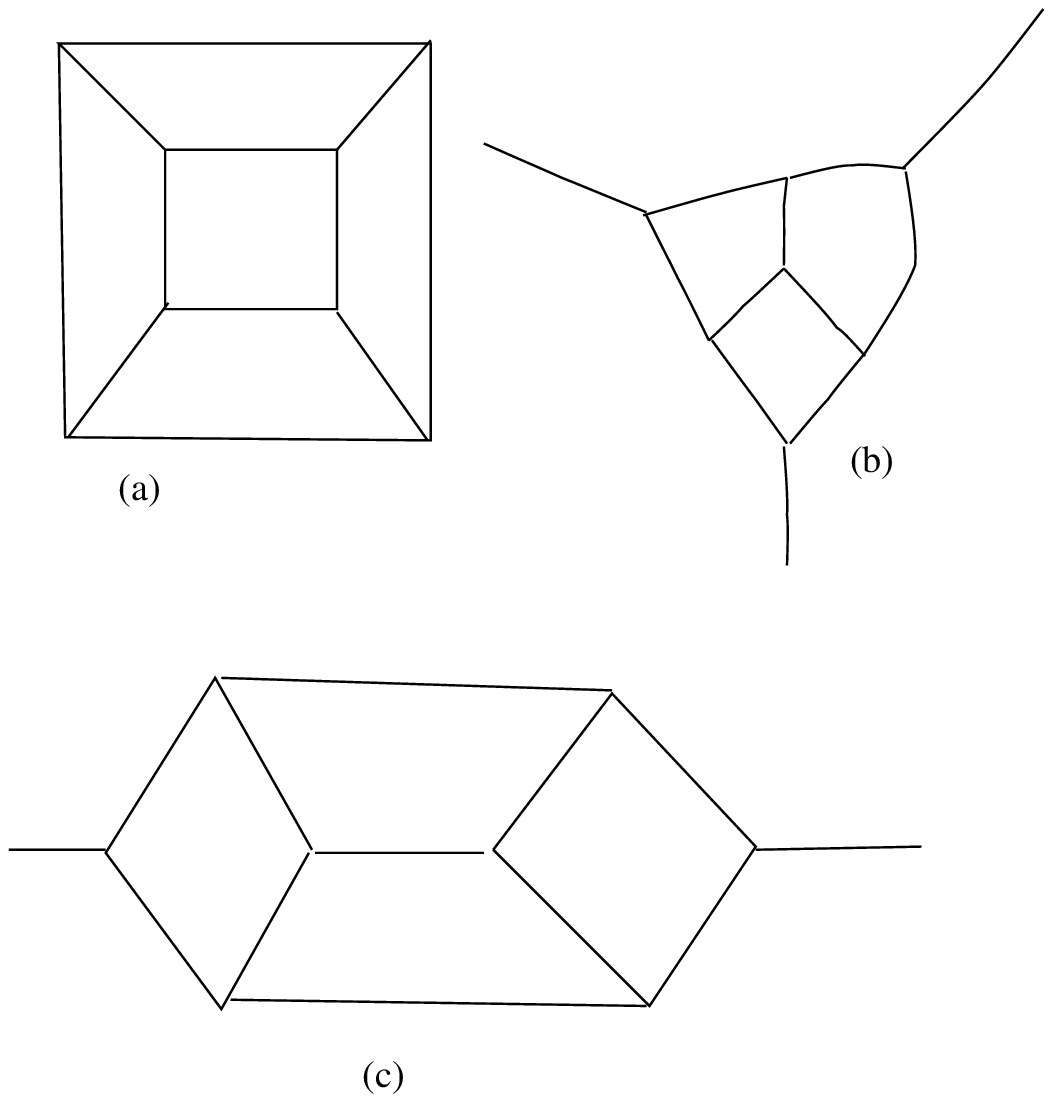}{three maps of the cube}{These three diagrams are maps
of the cube, stretched out in the plane.  In (a), a point has been
removed from a face in order to stretch it out. In (b), a vertex has
been removed.  In (c), a point has been removed from an edge.}

Depending on which point you remove from the sphere,
you can get different maps in the plane.
For instance,
figure \ref{three maps of the cube}
shows three ways of representing the map depicting the edges and
vertices of the cube in the plane;
these three different pictures arise according to whether the point you
remove lies in the middle of a face,
lies on an edge,
or coincides with one of the vertices of the cube.

\subsection{Euler numbers}
For the regular polyhedra,
the \dt{Euler number} $V-E+F$ takes on the value 2,
where $V$ is the number of vertices, $E$ is the number of edges,
and $F$ is the number of faces.

The Euler number (pronounced `oiler number')
is also called the \dt{Euler characteristic},
and it is commonly denoted by the Greek letter $\chi$ (pronounced `kai',
to rhyme with `sky'):
\[
\chi=V-E+F
.
\]

\subsection{Discussion}

This exercise is designed to investigate the extent to which it is true
that the Euler number of a polyhedron is always equal to 2.
We also want you to gain some experience with representing polyhedra in the
plane using maps,
and with drawing dual maps.

We will be distributing examples of different polyhedra.
\begin{enumerate}
\item
For as many of the polyhedra as you can,
determine the values of $V$, $E$, $F$, and the Euler number $\chi$.
\item
When you are counting the vertices and so forth,
see if you can think of more than one way to count them,
so that you can check your answers.
Can you make use of symmetry to simplify counting?
\item
The number $\chi$ is frequently very small compared with $V$, $E$, and $F$,
Can you think of ways  to find the value of $\chi$ without
having to compute $V$, $E$, and $F$, by
`cancelling out' vertices or faces
with edges?
This gives another way to check your work.
\end{enumerate}

The \dt{dual} of a map is a map you get by putting a vertex in the
each face, connecting the neighboring faces by new edges which
cross the old edges, and removing all the old vertices and
edges.  To the extent feasible,
draw a map in the plane of the polyhedron,
draw (in a different color) the dual map,
and draw a net for the polyhedron as well.

\section{Notation for some knots}

It is a hard mathematical question to completely codify
all possible knots.  Given two knots, it is hard to tell whether
they are the same. It is harder still to tell for sure
that they are different.  

Many simple knots can be arranged in a certain form, as illustrated
below, which is described by a string of positive integers along
with a sign. 

\fig{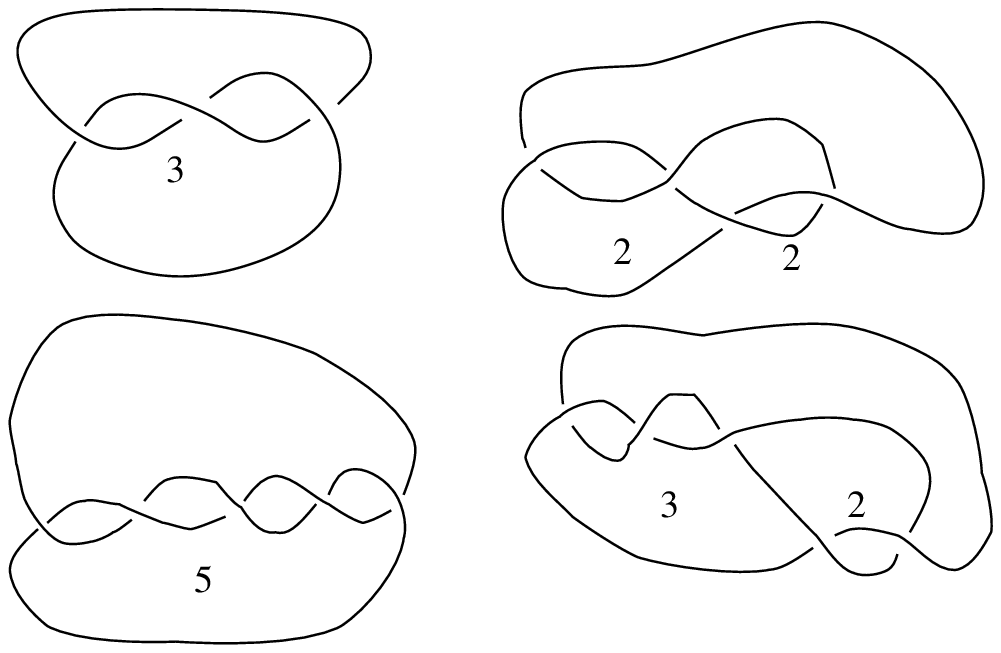}{Notation for certain knots}{Here are drawings of
some examples of knots that Conway `names' by a string of positive integers.
The drawings use the convention that when one strand crosses under
another strand, it is broken.
Notice that as you
run along the knot, the strand alternates going over and under at
its crossings.  Knots with this property are called alternating
knots.  Can you find any examples of knots with more than one
name of this type?}

\fig{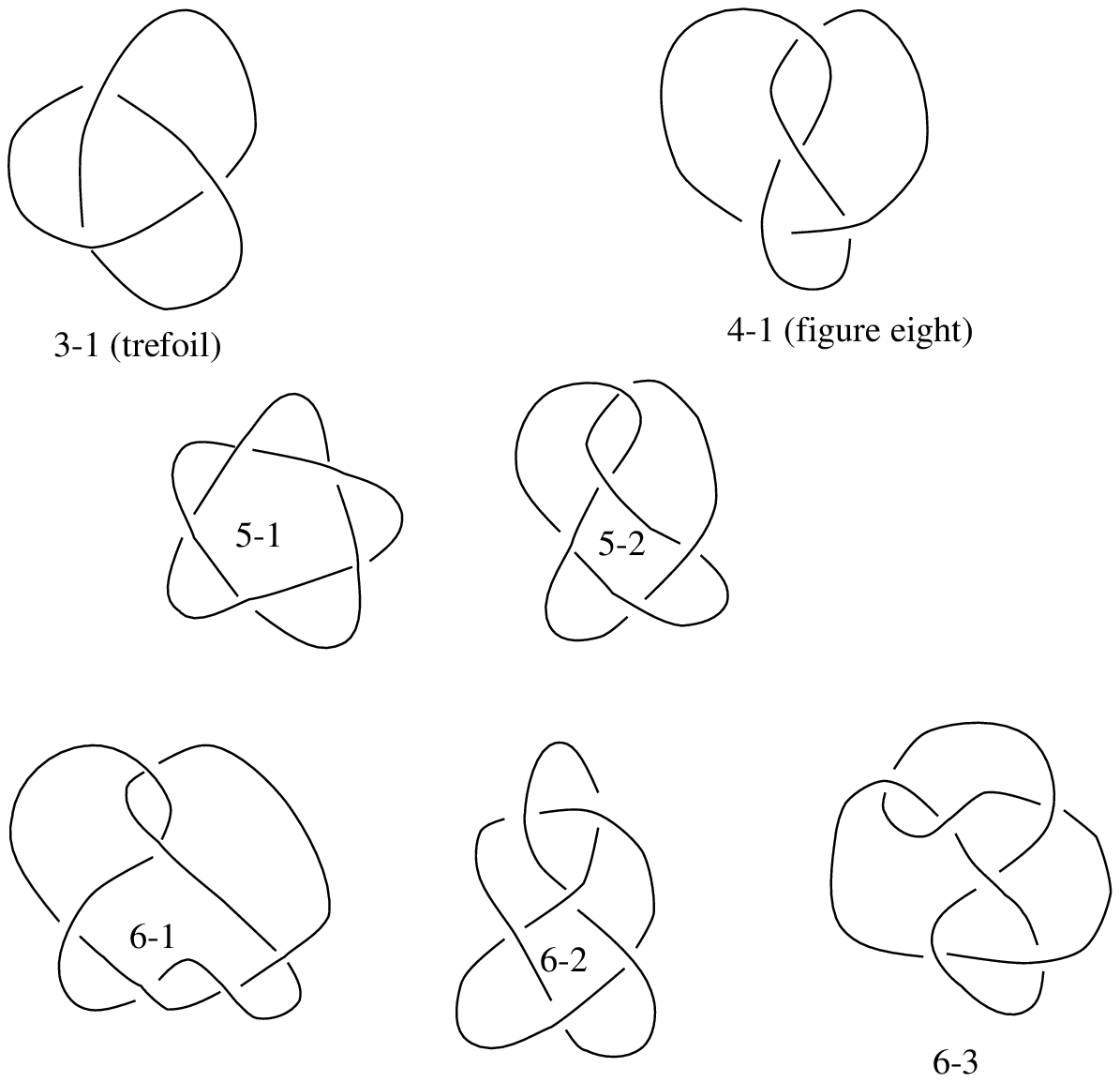}{knots up to six crossings}{Here are the knots with
up to six crossings.  The names follow an old system, used widely
in knot tables, where the $k$th knot with $n$ crossings is called
$n-k$.
  Mirror images are not included: some of these
knots are equivalent to their mirror images, and some are not.
Can you tell which are which?}

\section{Knots diagrams and maps}

A knot diagram gives a map on the plane, where there are four edges
coming together at each vertex.  Actually, it is better to think
of the diagram as a map on the sphere, with a polygon on the outside.
It sometimes helps in recognizing when diagrams are topologically
identical to label the regions with how many edges they have.

\section{Unicursal curves and knot diagrams}

A unicursal curve in the plane is a curve that you get when you put down
your pencil, and draw until you get back to the starting point.  As you
draw, your pencil mark can intersect itself, but you're not supposed to
have any triple intersections.  You could say that you pencil is allowed to
pass over an point of the plane at most twice.
This property of not having any triple intersections is \dt{generic}:
If you scribble the curve with your eyes closed (and somehow magically
manage to make the curve finish off exactly where it began),
the curve won't have any triple intersections.

A unicursal curve differs from the curves shown in knot diagrams in that
there is no sense of the curve's crossing over or under itself at an
intersection.  You can convert a unicursal curve into a knot diagram by
indicating (probably with the aid of an eraser), which strand crosses over and
which strand crosses under at each of the intersections.

A unicursal curve with 5 intersections can be converted into a knot diagram
in $2^5$ ways, because each intersection can be converted into a crossing
in two ways.
These 32 diagrams will not represent 32 different knots, however.

\subsection{Assignment}

\begin{enumerate}
\item
Draw the 32 knot diagrams that arise from the unicursal curve underlying
the diagram of knot 5-2 shown in the previous section, and identify
the knots that these diagrams represent.
\item
Show that any unicursal curve can be converted into a diagram of the
unknot.
\item
Show that any unicursal curve can be converted into the diagram of an
alternating knot in precisely two ways.
These two diagrams may or may knot represent different knots.
Give an example where the two knots are the same, and another where
the two knots are different.
\item
Show that any unicursal curve gives a map of the plane whose regions
can be colored black and white 
in such a way that adjacent regions have different colors.
In how many ways can this coloring be done?
Give examples.
\end{enumerate}

\section{Gas, water, electricity}

The diagram below shows three houses, each connected up to three utilities.
Show that it isn't possible to rearrange the connections so that they don't
intersect each other.  Could you do it if the earth were a not a sphere but
some other surface?

\fig{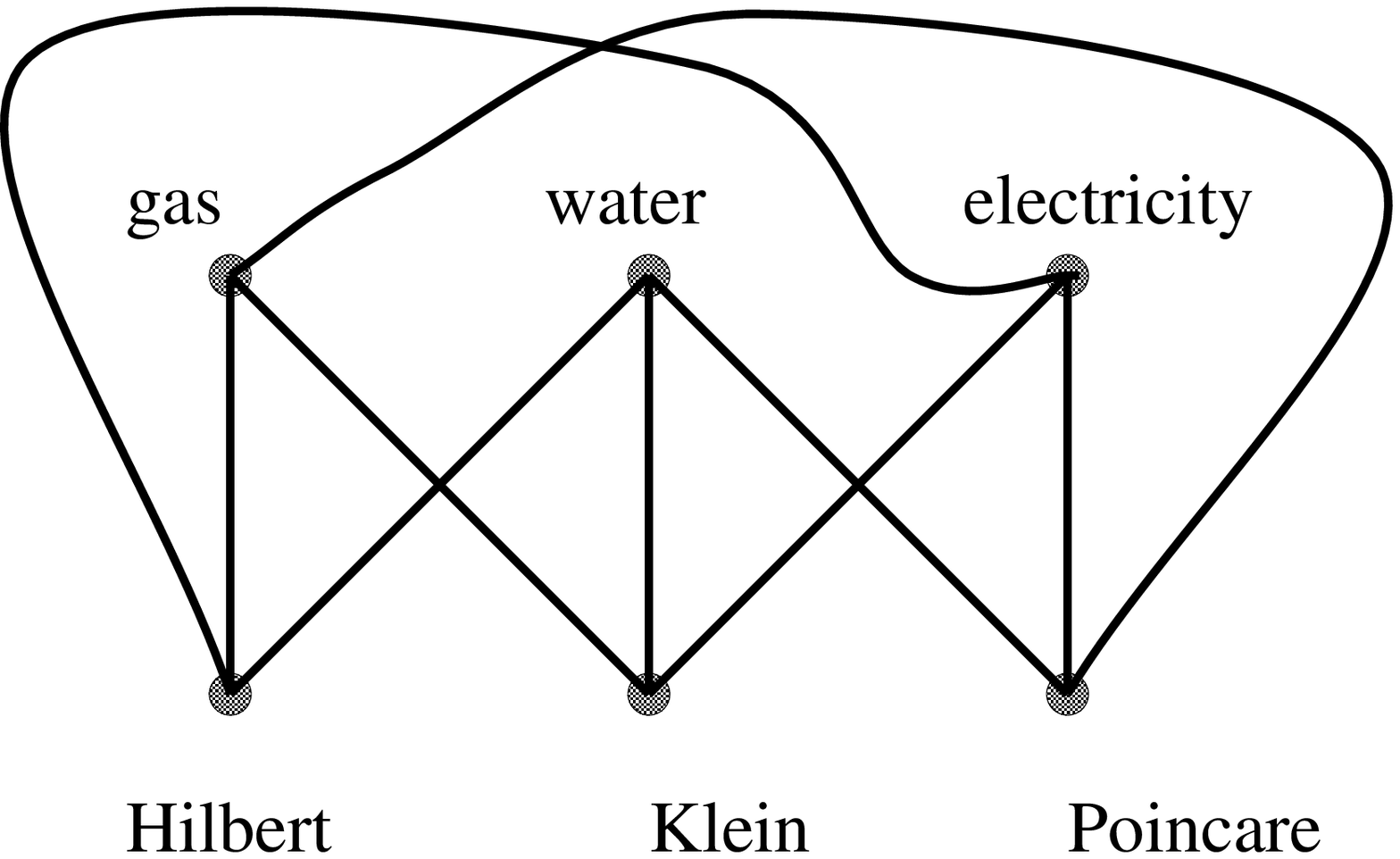}{Gas, water, electricity}{This is no good because we don't want
the lines to intersect.}

\section{Topology}

Topology is the theory of shapes which are allowed to stretch, compress,
flex and bend, but without tearing or gluing.  For example, a square
is topologically equivalent to a circle, since a square can be continously
deformed into a circle.  As another example, a doughnut and a coffee cup
with a handle for are topologically equivalent, since a doughnut can
be reshaped into a coffee cup without tearing or gluing.

\subsection{Letters}
As a starting exercise in topology, let's look at the letters of the
alphabet.   We think of the letters as figures made from 
lines and curves, without fancy doodads such as serifs.

{\bf Question.}
Which of the capital letters are topologically the same, and which
are topologically different?
How many topologically different capital letters are there?

\section{Surfaces}

A \dt{surface}, or \dt{2-manifold},
is a shape any small enough neighborhood of which
is topologically equivalent to a neighborhood of a point in the plane.
For instance, a the surface of a cube is a surface topologically
equivalent to the surface of a sphere.
On the other hand, if we put an extra wall inside a cube dividing it
into two rooms, we no longer have a surface, because there are points
at which three sheets come together. No small neighborhood of those
points is topologically equivalent to a small neighborhood in the plane.

\fig{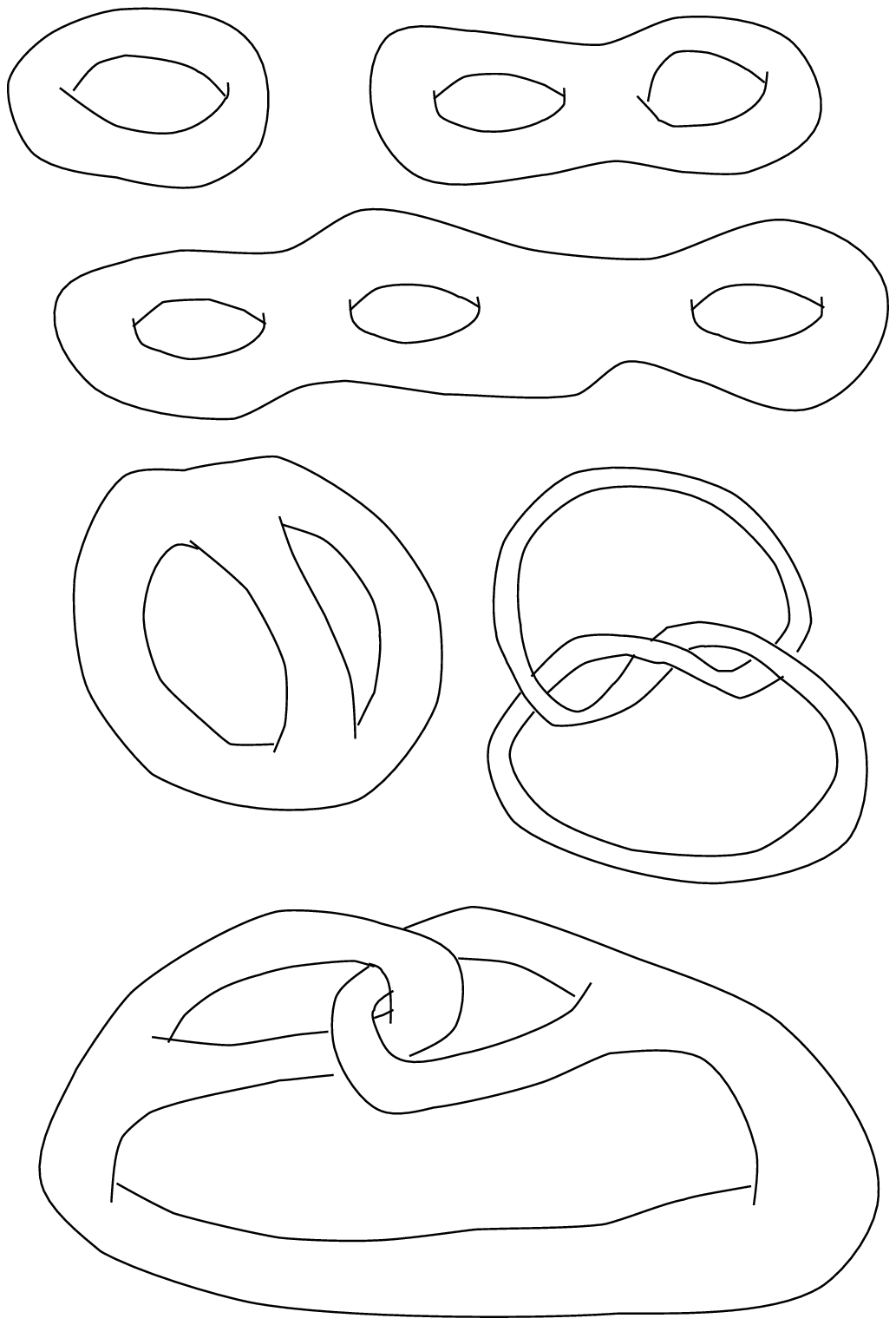}{Some surfaces}{Here are some pictures of surfaces.  The
pictures are intended to indicate things like doughnuts and pretzels rather than flat strips of paper.
Can you identify these surfaces, topologically?  Which ones are topologically
the same intrinsically, and which extrinsically?}
Recall that you get a torus by identifying the sides of a rectangle as
in Figure 2.10 of \em{SS} \em{(The Shape of Space)}.
If you identify the sides slightly differently,
as in Figure 4.3, you get a surface called a \dt{Klein bottle},
shown in Figure 4.9.

\subsection{Discussion}

\begin{enumerate}
\item
Take some strips and join the opposite ends of each strip together as follows:
with no twists;
with one twist (half-turn)---this is called a \dt{M\"obius strip};
with two twists;
with three twists.
\item
Imagine that you are a two-dimensional being who lives in one of these four
surfaces.
To what extent can you tell exactly which one it is?
\item
Now cut each of the above along the midline of the original strip.
Describe what you get. Can you explain why?
\item
What is the Euler number of a disk?  A M\"obius strip?
A torus with a circular hole cut from it?
A Klein bottle?
A Klein bottle with a circular hole cut from it?
\item
What is the maximum number of points in the plane such that you can draw
non-intersecting segments joining each pair of points?
What about on a sphere?  On a torus?
\end{enumerate}

\section{How to knit a M\"obius Band}

Start with a different color from the one you want to make the band in.  Call
this the spare color.  With the spare color and normal knitting needles
cast on 90 stitches.

Change to your main color yarn.  Knit your row of 90 stitches onto a circular
needle.  Your work now lies on about 2/3 of the needle.  One end of the work
is near the tip of the needle and has the yarn attached.  This is the
working end.  Bend the working end around to the other end of your work,
and begin to knit those stitches onto the working end, but 
{\em do not \/} slip them off the other end of the needle as you
normally would.  When you have knitted all 90 stitches in this way, the
needle loops the work twice.

Carry on knitting in the same direction, slipping stitches off the needle
when you knit them, as normal.  The needle will remain looped around the 
work twice.  Knit five `rows' (that is $5 \times 90$ stitches) in this way.

Cast off.  You now have a Mobius band with a row of your spare color
running around the middle.  Cut out and remove the spare colored yarn.
You will be left with one loose stitch in
your main color which needs to be secured.

\fig{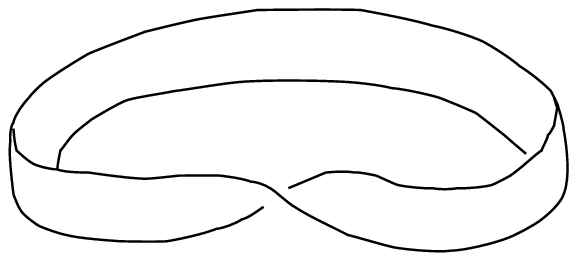}{Mobius band}{A Mobius band.}

(Expanded by Maria Iano-Fletcher from an original recipe by Miles Reid.)

\section{Geometry on the sphere}

We want to explore some aspects of geometry on the surface of the sphere.
This is an interesting subject in itself,
and it will come in handy later on when we discuss Descartes's angle-defect
formula.

\subsection{Discussion}

Great circles on the sphere are the analogs of straight lines in the plane.
Such curves are often called \dt{geodesics}.  A \dt{spherical triangle} is
a region of the sphere bounded by three arcs of geodesics.

\begin{enumerate}
\item
Do any two distinct points on the sphere determine a unique geodesic?
Do two distinct geodesics intersect in at most one point?
\item
Do any three `non-collinear' points on the sphere determine a unique
triangle?
Does the sum of the angles of a spherical triangle always equal $\pi$?
Well, no.  What values can the sum of the angles take on?
\end{enumerate}

\vskip 1in
The area of a spherical triangle is the amount by which the sum of its
angles exceeds the sum of the angles ($\pi$) of a Euclidean triangle.
In fact, for any spherical polygon, the sum of its angles minus the
sum of the angles of a Euclidean polygon with the same number of sides
is equal to its area.  

A proof of the area formula can be found in Chapter 9 of Weeks, {\it The
Shape of Space}.

\section{Course projects}

We expect everyone to do a project for the course.  On the last day of
the course, Friday, June 28th, we will hold a Geometry Fair, where
projects will be exhibited.  Parents and any other interested people
are invited.  

Here are some ideas, to get you started thinking about possible projects.
Be creative---don't feel limited by these ideas.

\begin{itemize}
\item  Write a computer program that allows the user to
select one of the $17$ planar symmetry groups, start doodling,
and see the pattern replicate, as in Escher's drawings.

\item  Write a similar program for drawing tilings of the
hyperbolic plane, using one or two of the possible hyperbolic symmetry groups.

\item  Make sets of tiles which exhibit various kinds of
symmetry and which tile the plane in various symmetrical patterns.

\item  Write a computer program that replicates three-dimensional
objects according to a three-dimensional pattern, as in the tetrahedron,
octahedron, and icosahedron.

\item  Construct kaleidoscopes for tetrahedral, octahedral and
icosahedral symmetry.

\item Construct a four-mirror kaleidoscope, giving a three-dimensional
pattern of repeating symmetry.

\item
The Archimidean solids are solids whose faces are regular
polygons (but not necessarily all the same) such that every vertex is symmetric
with every other vertex.  Make models of the the Archimedean solids

\item Write a computer program for visualizing four-dimensional
space.

\item Make stick models of the regular four-dimensional solids.

\item Make models of three-dimensional cross-sections of 
regular four-dimensional solids.

\item Design and implement three-dimensional tetris.

\item Make models of the regular star polyhedra
(Kepler-Poinsot polyhedron).

\item
Knit a Klein bottle, or a projective plane.

\item
Make some hyperbolic cloth.

\item
Sew topological surfaces and maps.

\item Infinite Euclidean polyhedra.

\item Hyperbolic polyhedra.

\item  Make a (possibly computational) orrery.

\item  Design and make a sundial.

\item  Astrolabe  (Like a primitive sextant).

\item Calendars: perpetual, lunar, eclipse.

\item Cubic surface with 27 lines.

\item Spherical Trigonometry or Geometry: Explore spherical trigonometry or
geometry. What is the analog
on the sphere of a circle in the plane? Does every spherical triangle
have a unique inscribed and circumscribed circle? Answer these and other similar questions.

\item Hyperbolic  Trigonometry or Geometry: Explore hyperbolic trigonometry or
 geometry. What is the analog
in the hyperbolic plane  of a circle in the Euclidean plane? Does every hyperbolic  triangle
have a unique inscribed and circumscribed circle? Answer these and other similar questions.

\item Make a convincing model showing how a torus can be filled
with circular circles in four different ways.

\item Turning the sphere inside out.

\item  Stereographic lamp.

\item  Flexible polyhedra.

\item  Models of ruled surfaces.

\item   Models of the projective plane.

\item  Puzzles and models illustrating extrinsic topology.

\item Folding ellipsoids, hyperboloids, and other figures.

\item Optical models: elliptical mirrors, {\it etc.}

\item Mechanical devices for angle trisection, {\it etc.}

\item Panoramic polyhedron (similar to an astronomical globe)
made from faces which are photographs.
\end{itemize}

\section{The angle defect of a polyhedron}

The \dt{angle defect} at a vertex of a polygon is defined to be $2 \pi$
minus the sum of the angles at the corners of the faces at that vertex.
For instance, at any vertex of a cube there are three angles of $\pi/2$,
so the angle defect is $\pi/2$.  You can visualize the angle defect by
cutting along an edge at that vertex, and then flattening out a neighborhood
of the vertex into the plane.  A little gap will form where the slit is:
the angle by which it opens up is the angle defect.

The \dt{total angle defect} of the polyhedron is gotten by adding up the
angle defects at all the vertices of the polyhedron. For a cube,
the total angle defect is $8 \times \pi/2 = 4 \pi$.
\subsection{Discussion}
\begin{enumerate}
\item
What is the angle sum for a polygon (in the plane) with $n$ sides?
\item
Determine the total angle defect for each of the 5 regular polyhedra,
and for the polyhedra handed out.
\end{enumerate}

\section{Descartes's Formula.}

The \dt{angle defect} at a vertex of a polygon was  defined to be the amount
by which the sum of the angles at the corners of the faces at that vertex
falls short of $2\pi$  
and the \dt{total angle defect} of the polyhedron was defined to be what one got when one added up the
angle defects at all the vertices of the polyhedron.
We call the total defect $T$.
Descartes discovered that there is a connection between the total defect, $T$, and the Euler Number $E - V - F$. Namely,
\begin{equation}
T = 2\pi(V- E + F).
\end{equation}
Here are two proofs.  They both use the fact that the sum of the
angles of a polygon with $n$ sides is $(n-2)\pi$.

\subsection{First proof}
Think of $2 \pi (V - E + F)$ as putting $+ 2 \pi$ at each vertex,
$-2\pi$ on each edge, and $+2 \pi$ on each face.

We will try to cancel out the terms as much as possible, by grouping
within polygons.

For each edge, there is $-2\pi$ to allocate.  An edge has a polygon on
each side: put $-\pi$ on one side, and $-\pi$ on the other.

For each vertex, there is $+2\pi$ to allocate: we will do it according to
the angles of polygons at that vertex.  If the angle of a polygon at
the vertex is $a$, allocate $a$ of the $2 \pi$ to that polygon.
This leaves something at the vertex: the angle defect.

In each polygon, we now have a total of the sum of its angles minus
$n \pi$ (where $n$ is the number of sides) plus $2 \pi$. Since
the sum of the angles of any polygon is $(n-2)\pi$, this is 0.  
Therefore, 
$$
2 \pi (V-E+F) = T .
$$

\vskip 3in
\subsection{Second  proof}

We begin to compute:

$$T = \sum_{\hbox {Vertices}} \hbox{the angle defect at the vertex}.$$
$$ = \sum_{\hbox {Vertices}} (2\pi - \hbox{the sum of the angles at the corners of those faces that meet at the vertex}).$$
$$ = 2\pi V - \sum_{\hbox {Vertices}} (\hbox{the sum of the angles at the corners of those faces that meet at the vertex}).$$
$$ = 2\pi V - \sum_{\hbox {Faces}} \hbox {the sum of the interior angles of the face}.$$
$$ = 2\pi V - \sum_{\hbox{Faces}} (n_f - 2)\pi.$$
 Here $n_f$ denotes the number of edges on the face $f$. 
$$ T = 2\pi V - \sum_{\hbox{Faces}} n_f\pi + \sum_{\hbox{Each face}} 2\pi. $$ 
Thus 
$$ T = 2\pi V -( \sum_{\hbox{Faces}} \hbox{the number of edges on the face}\cdot\pi) + 2\pi F.$$
If we sum  the number of edges  on each face over all of the faces, we will have counted each edge {\bf twice}. Thus
$$ T = 2\pi V - 2E\pi + 2\pi F.$$ Whence,
$$ T = 2\pi(V - E + F).$$

\subsection{Discussion}
Listen to both proofs  given in class.
\begin{enumerate}
\item Discuss both proofs with the aim of understanding them.
\item Draw a sketch of the first proof in the blank space above.
\item Discuss the differences between the two proofs. Can you describe the
ways in which they are different? Which of you feel 
the first is easier to understand?  Which of you feel the second is easier 
to understand?  Which is more pleasing? Which is more conceptual? 
\end{enumerate}

\section{Exercises in imagining}

How do you imagine geometric figures in your head?  Most people talk about
their three-dimensional imagination as `visualization', but that isn't
exactly right.  A visual image is a kind of picture, and it is really
two-dimensional.  The image you form in your head is more conceptual
than a picture---you locate things in more of a three-dimensional model than
in a picture.  In fact, it is quite hard to go from a mental image to a
two-dimensional visual picture.  Children struggle
long and hard to learn to draw because of the real conceptual difficulty
of translating three-dimensional mental images into
two-dimensional images.

Three-dimensional mental images are connected with your visual sense,
but they are also connected with your sense of place and motion. 
In forming an image, it often helps to imagine
moving around it, or tracing it out with your hands.  The size of an
image is important.  Imagine a little half-inch sugarcube in your hand,
a two-foot cubical box, and a ten-foot cubical room that you're inside.
Logically, the three cubes have the same information,
but people often find it easier to manipulate the larger image that they
can move around in.

Geometric imagery is not just something that you are either born with or
you are not.  Like any other skill, it develops with practice. 

Below are some images to practice with.  Some are two-dimensional, some are
three-dimensional.  Some are easy, some are hard, but not necessarily in
numerical order.  Find another person to work with in going through
these images.  Evoke the images by talking about them, not by drawing them.
It will probably help to close your
eyes, although sometimes gestures and drawings in the air will help.
Skip around to try to find exercises that are the right level for you.

When you have gone through these images and are hungry for more,
make some up yourself.

\begin{enumerate}

\item Picture your first name, and read off the letters
backwards.  If you can't see your whole name at once,
do it by groups of three letters.  Try
the same for your partner's name, and for a few other words.
Make sure to do it by sight, not by sound.

\item
Cut off each corner of a square, as far as the midpoints of the edges.
What shape is left over?  How can you re-assemble the four corners to make 
another square?

\item
Mark the sides of an equilateral triangle into thirds.  Cut off each corner
of the triangle, as far as the marks.  What do you get?

\item
Take two squares.  Place the second square centered over the first square
but at a forty-five degree angle.  What is the intersection of the two
squares?

\item
Mark the sides of a square into thirds, and cut off each of its corners
back to the marks.  What does it look like?

\item
How many edges does a cube have?

\item
Take a wire frame which forms the edges of a cube.  Trace out
a closed path which goes exactly once through each corner.

\item
Take a $3 \times 4$ rectangular array of dots in the plane,
and connect the dots vertically and horizontally.  How many squares
are enclosed?

\item
Find a closed
path along the edges of the diagram above which visits each
vertex exactly once?  Can you do it for a $3 \times 3$ array of dots?

\item
How many different colors are required to color the faces of a cube so that
no two adjacent faces have the same color?

\item
A tetrahedron is a pyramid with a triangular base.  How many faces
does it have?  How many edges?  How many vertices?

\item
Rest a tetrahedron on its base, and cut it halfway up.  What shape is the
smaller piece?   What shapes are the faces of the larger pieces?

\item Rest a tetrahedron so that it is balanced on one edge,
and slice it horizontally halfway between its lowest edge and its highest
edge.  What shape is the slice?

\item
Cut off the corners of an equilateral triangle as far as the midpoints of its
edges.  What is left over?

\item
Cut off the corners of a tetrahedron as far as the midpoints of the edges.
What shape is left over?

\item
You see the silhouette of a cube, viewed from the corner.  What does
it look like?

\item
How many colors are required to color the faces of an octahedron so that
faces which share an edge have different colors?

\item
Imagine a wire is shaped to go up one inch, right one inch, back one inch,
up one inch, right one inch, back one inch, \ldots .   What does it look
like, viewed from different perspectives?

\item  The game of tetris has pieces whose shapes are all the
possible ways that four squares can be glued together along edges.
Left-handed and right-handed forms are distinguished.  What are the shapes,
and how many are there?

\item  Someone is designing a three-dimensional tetris, and wants
to use all possible shapes formed by gluing four cubes together.  What
are the shapes, and how many are there?

\item  An octahedron is the shape formed by gluing together
equilateral triangles four to a vertex.  Balance it on a corner, and
slice it halfway up.  What shape is the slice?

\item  Rest an octahedron on a face, so that another face is
on top.  Slice it halfway up.  What shape is the slice?

\item Take a $3 \times 3 \times 3$ array of dots in space,
and connect them by edges up-and-down, left-and-right, and forward-and-back.
Can you find a closed path which visits every dot but one exactly once?
Every dot?

\item  Do the same for a $4 \times 4 \times 4$ array of dots,
finding a closed path that visits every dot exactly once.

\item
What three-dimensional solid has circular profile viewed from above,
a square profile viewed from the front, and a triangular profile viewed
from the side?   Do these three profiles determine the three-dimensional
shape?

\item
Find a path through edges of the dodecahedron which visits each vertex
exactly once.
\end{enumerate}

\section{Curvature of surfaces}

If you take a flat piece of paper and bend it gently, it bends in only one
direction at a time.  At any point on the paper, you can find 
at least one direction through which there is a straight line on the surface.
You can bend it into a cylinder, or into a cone, but 
you can never bend it without crumpling or distorting to the get
a portion of the surface of a sphere. 

If you take the skin of a sphere, it cannot be flattened out into the plane
without distortion or crumpling.  This phenomenon is familiar
from orange peels or apple peels.  Not even a small area of the skin
of a sphere can be flattened out without some distortion, although the
distortion is very small for a small piece of the sphere.  That's why
rectangular maps of small areas of the earth work pretty well, but
maps of larger areas are forced to have considerable distortion.

The physical descriptions of what happens as you bend various surfaces without
distortion do not have to do with the topological properties of the surfaces.
Rather, they have to do with the \dt{intrinsic geometry} of the surfaces.
The intrinsic geometry has to do with geometric properties which can
be detected by measurements along the surface, without considering
the space around it.

\bigskip
There is a mathematical way to explain the intrinsic geometric property
of a surface that tells when one surface can or cannot be bent into another.
The mathematical concept is called the \dt{Gaussian curvature} of a surface,
or often simply the \dt{curvature} of a surface.
This kind of curvature is not to be confused with the curvature
of a curve.  The curvature of a curve is an extrinsic geometric property,
telling how it is bent in the plane, or bent in space.
Gaussian curvature is an intrinsic geometric
property: it stays the same no matter how a surface is bent, as long
as it is not distorted, neither stretched or compressed.

To get a first qualitative idea of how curvature works, here are some
examples.

A surface which bulges out in all directions, such as the surface
of a sphere, is \dt{positively curved}.  A rough test for positive curvature
is that if you take any point on the surface, there is some plane touching
the surface at that point so that the surface lies all on one
side except at that point.  No matter how you (gently) bend the surface,
that property remains.

A flat piece of paper, or the surface of a cylinder or cone, has 0 curvature.

  A saddle-shaped
surface has negative curvature: every plane through a point on the saddle
actually cuts the saddle surface in two or more pieces.

{\bf Question.}
What surfaces can you think of that have positive,
zero, or negative curvature.

Gaussian curvature is a numerical quantity associated to an
area of a surface, very closely related to angle defect.  Recall
that the angle defect of a polyhedron at a vertex
is the angle by which a small neighborhood of a vertex opens up,
when it is slit along one of the edges going into the vertex.

The total Gaussian curvature of a region on a surface is the angle by which
its boundary opens up, when laid out in the plane.  To actually measure 
Gaussian curvature of a region bounded by a curve,
you can cut out a narrow strip on the surface in neighborhood of the
bounding curve.  You also need to cut open the curve, so it will be
free to flatten out.  Apply it to a flat surface, being careful to
distort it as little as possible.  If the surface is positively curved
in the region inside the curve, when you flatten it out, the curve will
open up.  The angle between the tangents to the curve at the two sides
of the cut is the total Gaussian curvature.  This is like angle defect:
in fact, the total curvature of a region of a polyhedron containing exactly
one vertex is the angle defect at that vertex.  You must pay attention
pay attention not just to the angle between the ends of the strip, but
how the strip curled around,
keeping in mind that the standard for zero curvature is
a strip which comes back and meets itself.  Pay attention to $\pi$'s and
$2\pi$'s.

\fig{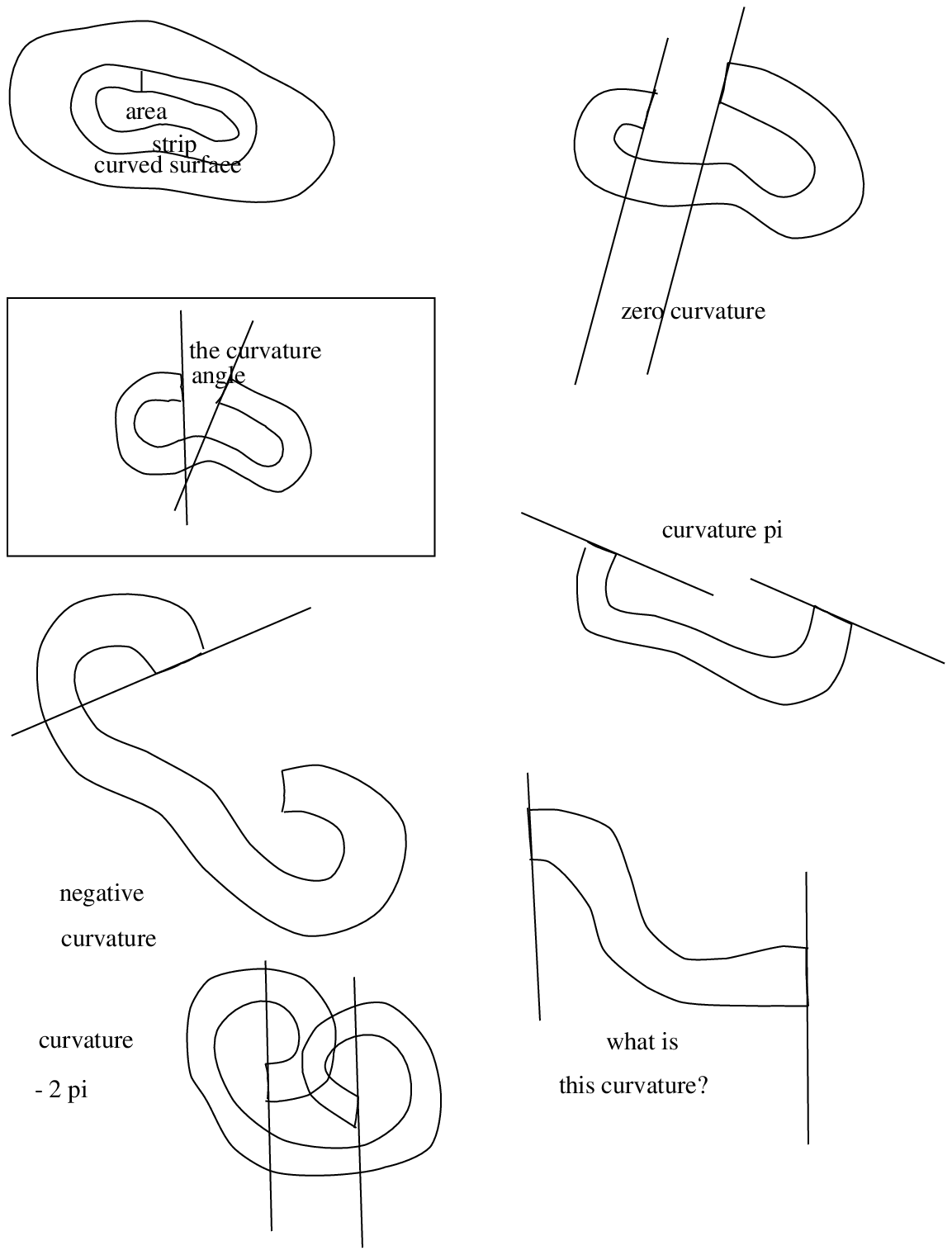}{measuring curvature}{This diagram illustrates how to measure
the total Gaussian curvature of a patch by cutting out a strip
which bounds the patch, and laying it out on a flat surface.  The angle
by which the strip `opens up' is the total Gaussian curvature.  You must
pay attention not just to the angle between the lines on the paper, but
how it got there, keeping in mind that the standard for zero curvature is
a strip which comes back and meets itself.  Pay attention to $\pi$'s and
$2\pi$'s.}

If the total curvature inside the region is negative, the strip will
curl around further than necessary to close.  The curvature is negative,
and is measured by the angle by which the curve overshoots.

A less destructive way to measure total Gaussian curvature of a region
is to apply narrow strips of paper to the surface, e.g., masking tape.
They can be then be removed and flattened out in the plane to measure
the curvature.

{\bf Question.}
Measure the total Gaussian curvature of 
\begin{enumerate}
\item
a cabbage leaf.
\item
a kale leaf
\item
a piece of banana peel
\item
a piece of potato skin
\end{enumerate}
If you take two adjacent regions, is
the total curvature in the whole equal to the sum of the total curvature
in the parts?  Why?

The angle defect of a convex polyhedron at one of its vertices can be
measured by rolling the polyhedron in a circle around its vertex.
Mark one of the edges, and rest it on a sheet of paper.  Mark the line
on which it contacts the paper.  Now roll the polyhedron, keeping 
the vertex in contact with the paper.  When the given edge first touches
the paper again, draw another line.  The angle between the two lines
(in the area where the polyhedron did not touch) is the angle defect.
In fact, the area where the polyhedron did touch the paper can be rolled
up to form a paper model of a neighborhood of the vertex in question.

A polyhedron can also be rolled in a more general way.  Mark some
closed path on the surface of the polyhedron, avoiding vertices.
Lay the polyhedron on a sheet of paper so that part of the curve
is in contact.  Mark the position of one of the edges in contact with
the paper.
now roll the polyhedron, along the curve, until the original face
is in contact again, and mark the new position of the same edge.
What is the angle between the original position of the line, and the
new position of the line?

\section{Gaussian curvature}
\subsection{Discussion}
\begin{enumerate}
\item
What is the curvature inside the region on a sphere {\em exterior}
to a tiny circle?
\item
On a polyhedron,
what is the curvature inside a region containing a single vertex?
two vertices?
all but one vertex?
all the vertices?
\end{enumerate}

\section{Clocks and curvature}

The total curvature of any surface topologically equivalent to the sphere
is $4\pi$.
This can be seen very simply
from the definition of the curvature of a region
in terms of the angle of rotation when the surface
is rolled around on the plane;
the only problem is the perennial one of
keeping proper track of multiples of $\pi$ when measuring the angle of
rotation.
Since are trying to show that the total curvature is a specific multiple of
$\pi$,
this problem is crucial.
So to begin with let's think carefully about
how to reckon these angles correctly
\subsection{Clocks}
Suppose we have a number of clocks on the wall.
These clocks are good mathematician's clocks, with a 0 up at the top where
the 12 usually is.
(If you think about it, 
0 o'clock makes a lot more sense than 12 o'clock:
With the 12 o'clock system,
a half hour into the new millennium on 1 Jan 2001,
the time will be 12:30 AM, the 12 being some kind of hold-over from the
departed millennium.)

Let the clocks be labelled $A$, $B$, $C$, \ldots.
To start off, we set all the clocks to 0 o'clock.
(little hand on the 0; big hand on the 0),
Now we set clock $B$ ahead half an hour
so that it now the time it tells is 0:30
(little hand on the 0 (as they say); big hand on the 6).
What angle does its big hand make with that of clock $A$?
Or rather, through what angle has its big hand moved relative to that
of clock $A$?
The angle is $\pi$.
If instead of degrees or radians, we measure our angles in \dt{revs}
(short for \dt{revolutions}),
then the angle is $1/2$ rev.
We could also say that the angle is $1/2$ hour:
as far as the big hand of a clock is concerned,
an hour is the same as a rev.

Now take clock $C$ and set it to 1:00.
Relative to the big hand of clock $A$, the big hand of $C$ has moved
through an angle of $2 \pi$, or 1 rev, or 1 hour.
Relative to the big hand of $B$, the big hand of $C$ has moved through
an angle of $\pi$, or $1/2$ rev.
Relative to the big hand of $C$, the big hand of $A$ has moved through
an angle of $-2\pi$, or $-1$ rev,
and the big hand of $B$ has moved $- \pi$, or $-1$ rev.

\subsection{Curvature}
Now let's describe how to find the curvature inside a disk-like region $R$
on a surface $S$, i.e. a region topologically equivalent to a disk.
What we do is cut a small circular band running around the boundary of the
region, cut the band  open to form a thin strip, lay the thin strip 
flat on the plane,
and measure the angle between the lines at the two end of the strip.
In order to keep the $\pi$'s straight,
let us go through this process very slowly and carefully.

To begin with, let's designate the two ends of the strip
as the {\em left end} and the {\em right end} in such a way that
traversing the strip from the left end to the right end
corresponds to circling {\em clockwise} around the region.
We begin by fixing the left-hand end of the strip to the wall so that
the straight edge of the cut at the left end of the strip---the cut that
we made to convert the band into a strip---runs
straight up and down, parallel to the big hand of clock $A$,
and so that the strip runs off toward the right.
Now we move from left to right along the strip,
i.e. clockwise around the boundary of the region,
fixing the strip so that
it lies as flat as possible, until we come to the right end of the strip.
Then we look at the cut bounding the right-hand end of the strip,
and see how far it has turned relative to the left-hand end of the strip.
Since we were so careful in laying out the left-hand end of the strip,
our task in reckoning the angle of the right-hand end of the strip
amounts to deciding what time you get if you think of the right-hand end
of the strip as the big hand of a clock.
The curvature inside the region will correspond to the amount by which the time
told by the right-hand end of the strip falls short of 1:00.

For instance, say the region $R$ is a tiny disk in the Euclidean plane.
When we cut a strip from its boundary and lay it out as described above,
the time told by its right hand end will be precisely 1:00, so the curvature
of $R$ will be exactly 0.
If $R$ is a tiny disk on the sphere, then when the strip is laid out the time
told will be just shy of 1:00, say 0:59, and the curvature of the region
will be $\frac{1}{60}$ rev, or $\frac{\pi}{30}$.

When the region $R$ is the lower hemisphere of a round sphere,
the strip you get will be laid out in a straight line,
and the time told by the right-hand end will be 0:00,
so the total curvature will be 1 rev, i.e. $2 \pi$.
The total curvature of the upper hemisphere is $2 \pi$ as well,
so that the total curvature of the sphere is $4 \pi$.

Another way to see that the total curvature of the sphere is $4 \pi$
is to take as the region $R$ the {\em outside} of a small circle on the
sphere.
When we lay out a strip following the prescription above, being sure to
traverse the boundary of the region $R$ in the clockwise sense as viewed
from the point of view of the region $R$, we see that the time told by the
right hand end of the strip is very nearly $-1$ o'clock!
The precise time will be just shy of this,
say $-1$:59,
and the total curvature of the region will then be $1 \frac{59}{60}$ revs.
Taking the limit,
the total curvature of the sphere is 2 revs, or $4 \pi$.

But this last argument will work equally well on any surface topologically
equivalent to a sphere, so any such surface has total curvature $4 \pi$.

\subsection{Where's the beef?}
This proof that the total curvature of a topological sphere is $4 \pi$ 
gives the definite feeling of being some sort of trick.
How can we get away without doing any work at all?
And why doesn't the argument work equally well on a torus, which as we know should have total curvature 0? What gives?

What gives is the lemma that states that if you take a disklike region $R$
and divide it into two disklike subregions $R_1$ and $R_2$, then the
curvature inside $R$ when measured by laying out its boundary is the sum
of the curvatures inside $R_1$ and $R_2$ measured in this way.
This lemma might seem like a tautology.
Why should there be anything to prove here?
How could it fail to be the case that the curvature inside the whole is
the sum of the curvatures inside the parts?
The answer is,
it could fail to be the case by virtue of our having given a faulty definition.
When we define the curvature inside a region, we have to make sure that the
quantity we're defining has the additivity property, or the definition is no
good.
Simply calling some quantity the curvature inside the region will not make
it have this additivity property.
For instance, what if we had defined the curvature inside a region to be
$4 \pi$, no matter what the region?
More to the point, what if in the definition of the curvature inside
a region we had forgotten the proviso that the region $R$ be disklike?
Think about it.

\section{Photographic polyhedron}

As you stand in one place and look around, up, and down,
there is a sphere's worth of directions you can look.  One
way to record what you see would be to construct a big sphere,
with the image painted on the inside surface.  To see the world
as viewed from the one place, you would stand on a platform in
the center of the sphere and look around.   We will call this
sphere the \dt{visual sphere}.  You can imagine a sphere, like a planetarium,
with projectors projecting a seamless image.  The image might be created
by a robotic camera device, with video cameras pointing in enough directions
to cover everything.

{\bf Question.}
What is the geometric relation of objects in space to their images
on the visual sphere?
\begin{enumerate}
\item
Show that the image of a line is an arc of a great circle.  If the line
is infinitely long, how long (in degrees) is the arc of the circle?
\item
Describe the image of several parallel lines.
\item
What is the image of a plane?
\end{enumerate}

Unfortunately, you can't order spherical prints from most photographic
shops.  Instead, you have to settle for flat prints.  Geometrically,
you can understand the relation of a flat print to the `ideal' print
on a spherical surface by constructing a plane tangent to the sphere
at a point corresponding to the center of the photograph.  You can
project the surface of the sphere outward to the plane, by following
straight lines from the center of the sphere to the surface of the sphere,
and then outward to the plane.  From this,
you can see that given size objects on the visual sphere do not always come
out the same size on a flat print.  The further they are from the center
of the photograph, the larger they are on the print.

Suppose we stand in one place, and take several photographs that overlap,
so as to construct a panorama.  If the camera is adjusted in
exactly the same way for the various photographs, and the prints are made
in exactly the same way, the photographs can be thought of as coming
from rectangles tangent to a copy of the visual sphere, of some size.
The exact radius of this sphere, the \dt{photograph sphere}
depends on the focal length of the camera lens, the size of prints, {\it etc.},
but it should be the same sphere for all the different prints.

If we try to just overlap
them on a table and glue them together, the images will not match up
quite right: objects on the edge of a print are larger than objects
in the middle of a print, so they can never be exactly aligned.

Instead, we should try to find the line where two prints would intersect
if they were arranged to be tangent to the sphere.  This line is
equidistant from the centers of the two prints.  You can find it by
approximately aligning the two prints on a flat surface, draw the line
between the centers of the prints, and constructing the perpendicular
bisector.  Cut along this line on one of the prints.  Now find the corresponding
line on the other print.  These two lines should match pretty closely.
This process can be repeated: now that the two prints have a better match,
the line segment between their centers can be constructed more accurately,
and the perpendicular bisector works better.  

If you perform this operation for a whole collection of photographs, you
can tape them together to form a polyhedron.  The polyhedron should be
circumscribed about a certain size sphere.  It can give an excellent
impression of a wide-angle view of the scene.  If the photographs cover
the full sphere, you can assemble them so that the prints are face-outwards.
This makes a globe, analogous to a star globe.  As you turn it around, you
see the scene in different directions.  If the photographs cover a fair
bit less than a full sphere, you can assemble them face inwards.  This
gives a better wide-angle view.

One way to do this is just to take enough photographs that you cover a
certain area of the visual sphere, match them up, cut them out, and tape
them together.  The polyhedron you get in this way will probably not
be very regular.

By choosing carefully the directions in which you take photographs, you
could make the photographic polyhedron have a regular, symmetric structure.
Using an ordinary lens, a photograph is not wide enough to fill the face
of any of the 5 regular polyhedra.  

An \dt{Archi\-medean polyhedron} is a polyhedron such that every face is a regular
polygon (but not necessarily all the same), and every vertex is symmetric
with every other vertex.  For instance, the soccer ball polyhedron, or
truncated icosahedron, is Arch\-imedean.

{\bf Question.}
Show that every Archi\-medean polyhedron is inscribed in a sphere.

The \dt{dual Archi\-medean polyhedra} are polyhedra which are dual to
Archi\-medean polyhedra.

{\bf Question.}
\begin{itemize}
\item
Show that each of the dual Archi\-medean polyhedra can be circumscribed
about a sphere.  
\item
Which polyhedra will work well to make a photographic polyhedron?
\end{itemize}

\section{Mirrors}

\subsection{Discussion}
\begin{enumerate}
\item
How do you hold two mirrors so as to get an integral number of images of
yourself?
Discuss the handedness of the images.
\item
Set up two mirrors so as to make perfect kaleidoscopic patterns.
How can you use them to make a snowflake?
\item
Fold and cut hearts out of paper.
Then make paper dolls.
Then honest snowflakes.
\item
Set up three or more mirrors so as to make perfect kaleidoscopic patterns.
Fold and cut such patterns out of paper.
\item
Why does a mirror reverse right and left rather than up and down?
\end{enumerate}

\section{More paper-cutting patterns}

Experiment with the constructions below. Put the best examples
into your journal, along with comments that describe and explain
what is going on.  Be careful to make your examples large enough
to illustrate clearly the symmetries that are present.
Also make sure that your cuts are interesting enough so that
extra symmetries do not creep in.
Concentrate on creating a collection of examples that will get across
clearly what is going on,
and include enough written commentary to make a connected narrative.

\begin{enumerate}
\item
{\bf Conical patterns.}
Many rotationally-symmetric designs, like the twin blades
of a food processor, cannot be made by folding and cutting.  However,
they can be formed by wrapping paper into a conical shape.

Fold a sheet of paper in half, and then unfold.  Cut along
the fold to the center of the paper.   Now wrap the paper into a conical
shape, so that the cut edge lines up with the uncut half of the fold.
Continue wrapping, so that the two cut edges line up and the
original sheet of paper wraps two full
turns around a cone.
Now cut out any pattern you like from the 
cone. Unwrap and lay it out flat.   The resulting pattern
should have two-fold rotational symmetry.

Try other examples of this technique, and also try experimenting with
rolling the paper more than twice around a cone.
\item
{\bf Cylindrical patterns.}  Similarly, it is possible to make
repeating designs on strips.   If you roll a strip of paper into a
cylindrical shape, cut it, and unroll it, you should get a repeating pattern
on the edge.  Try it.
\item
{\bf M\"{o}bius patterns.}  A M\"{o}bius band is formed by taking
a strip of paper, and joining one end to the other with a twist
so that the left edge of the strip continues to the right.

Make or round up a strip of paper which is long compared to its width
(perhaps made from ribbon, computer paper, adding-machine rolls, 
or formed by joining several shorter strips together end-to-end).  Coil
it around several times around in a M\"{o}bius band pattern.  Cut out
a pattern along the edge of the M\"{o}bius band, and unroll.
\item
{\bf Other patterns.}   Can you come up with any other creative ideas
for forming symmetrical patterns?
\end{enumerate}

\section{Summary}

In the past week we have discussed a number of different topics, many of
which seemed to be unrelated. 
When we began last week, we said that we would jump around from topic to topic during the first few days so that you would become familiar with a number
of different ideas and examples.
What we want to do today is to show you that there really is a method
 to our madness and  that there is  a connection between 
these seemingly diverse bits  of mathematics  and that the connection 
is one of the most deep and beautiful ones in mathematics. Virtually
 any property (visual or otherwise) that one naively chooses as a way 
 to describe (and quantify) 
a surface is related in a simple way to any other property
 one naively chooses and duly quantifies.    
 Here is a list of some of the things we touched upon last week.
\begin{itemize}
\item The Euler Number

\item Flashlights
\item Proofs  of the angular defect formula
\item Maps on surfaces
\item Area of a spherical triangle
\item Cabbage
\item Curvature
\item The Gauss map
\item Handle, holes, surfaces
\item Kale
\item Orientability
\end{itemize}

\section{The Euler Number}
If we have a polyhedron, we can compute its Euler number,  
$\chi = V - E + F$. In fact, we computed  Euler numbers
 {\it ad delectam}.
Why did we do this? One reason 
 is that they are easy to compute. But that is not obviously  a compelling
reason for doing anything in mathematics. The real reason is that it
 is an invariant of the surface (it does not depend upon what map one puts
on the surface) and because it is connected to a whole array of other
 properties  a surface might have that one might notice while trying to describe
it.
\subsection{ Descartes's Formula} 
One easy example of this is Descartes' formula. If one looks at
a polyhedral surface and makes a naive attempt to
describe it visually, one might try to describe how {\sl pointy}
 the surface is. A more sophisticated way to describe how {\sl pointy}
 a surface is at a vertex is to compute the angular defect at the vertex, that is
$$2\pi-(\hbox{ the sum of the angles of the faces meeting at the vertex }).$$
When we investigated how {\sl pointy} a polyhedron was,  summing over all
of the vertices to obtain the total angular defect $T$, we discovered that there was a direct connection
 between pointyness and Euler Number:
$$ T = 2\pi (V - E + F) .$$

\subsection{The Gauss Map (Flashlights)} 
Although projecting Conway's image onto the celestial sphere was fun, again it
was not in and of itself a mathematically valuable exercise. The point
was to get a feel for the Gauss map. The Gauss map is used to project
a surface onto the celestial sphere. For a polyhedron, we saw that, if one traced a path that
remained on a flat face, the Gauss image of that path was really a point.
We saw that if we traced a path that went around a vertex, the Gauss image
was a spherical polygon. If three edges  met at the given 
vertex, the Gauss image traced out a spherical triangle whose interior
could be thought of as the image of that vertex. Moreover, the angles of the triangle were
the supplements of the vertex angles. Using the formula for
the area of a spherical triangle, namely
 $$(\hbox{ the sum of the angles }) - \pi,$$ if the vertex angles were $\alpha, \beta, \hbox { and } \gamma$, the area
of the Gauss image of a path around the vertex would be
$$ \{(\pi - \alpha) + (\pi - \beta) + (\pi-\gamma ) \} - \pi = 2\pi - (\alpha + \beta + \gamma).$$ 
The right hand side of this formula is just the angular defect at the vertex.
Thus if we add up the areas of the images of path about all of the
vertices, we obtain the total defect $T$ of the original  surface.
Since no other parts of the image contribute to the area, we have shown that
$$\hbox {\bf the area of the Gauss image } =  T.$$
Exploiting the earlier connection, we can also say
$$\hbox {\bf the area of the Gauss image } =  2\pi(V-E+F).$$
This is known as the Gauss-Bonnet formula.
\subsection{Curvature (Kale and cabbage)}
Again, cutting up kale and cabbage was fun and the tape of Thurston and
Conway sticking  potato peel to the chalkboard will become a classic, but
there was a serious mathematical purpose behind it. If one looks at a surface
and wants to try to describe it visually, one might want to describe it by
 telling  how {\sl curly} it is. While the surface of a
 cylinder, for example, does not look visually as though it
 curves and bends very much, the surface of a trumpet does.
 Peeling a surface, that is, removing  a thin
strip from around a portion of the surface and then
seeing how much the angle between the ends of the strip opens up
 (or closes around)
as it is 
laid flat  quantifies the {\sl curviness} of
 the portion of the surface surrounded  by the strip. Mathematically,
 this is called the {\bf integrated curvature} of that portion of the surface.

 When we sum over portions that amount to the whole surface,
we get the {\bf total Gaussian curvature} of the surface. 

\subsubsection{Curvature for Polyhedra}
Lets apply these ideas to a polyhedron. In particular, we might
 consider a strip of polyhedron peel that just goes around one vertex
of a polyhedron.
Then we would find that the path opens up by an angle equal to the defect at that vertex, and so for such a path
$$ \hbox{ the total curvature enclosed }$$
$$= \hbox{ the defect at the enclosed vertex }$$
$$= \hbox { the area of the Gauss image }.$$
For a path that goes around several vertices the curvature is the sum of the
defects of all the surrounded vertices. Thus for a polyhedron, 
$$ K = \hbox {\bf Total Gaussian Curvature} = T.$$

\subsubsection{Curvature on surfaces}
To pass from a polyhedral surface to a smooth surface and to define curvature 
with mathematical precision, one needs to use integration in the definition
for $K$. But the conceptual idea is still the same. Any curved surface can be approximated by a polyhedral one with lots and
lots of vertices. The curvature of the surface within a path (a smooth
piece of peel) is then very nearly equal to the sum of the defects 
at all the encircled vertices. By a technical limiting argument that
involves integrals to give a precise meaning to curvature, $K$,
we find that {\bf for any surface}
$$ K = \hbox {\bf Total Gaussian Curvature} = T.$$

\subsection{Discussion}
\begin{itemize}
\item There are many other connections between these four 
concepts. Can you suggest any more?
(This is also a discussion question for the gang of four.)
\item 
The number of handles on a surface is another visual characteristic.
How does this relate to the total curvature?
\end{itemize}

\section{Symmetry and orbifolds}
\def\degree{^\circ}

Given a symmetric pattern, what happens when you identify equivalent
points?  It gives an object with interesting topological and geometrical
properties, called an orbifold.

The first instance of this is an object with bilateral symmetry, such
as a (stylized) heart.  Children learn to cut out a heart by folding a
sheet of paper in half, and cutting out half of the pattern.  When
you identify equivalent points, you get half a heart.

\fig{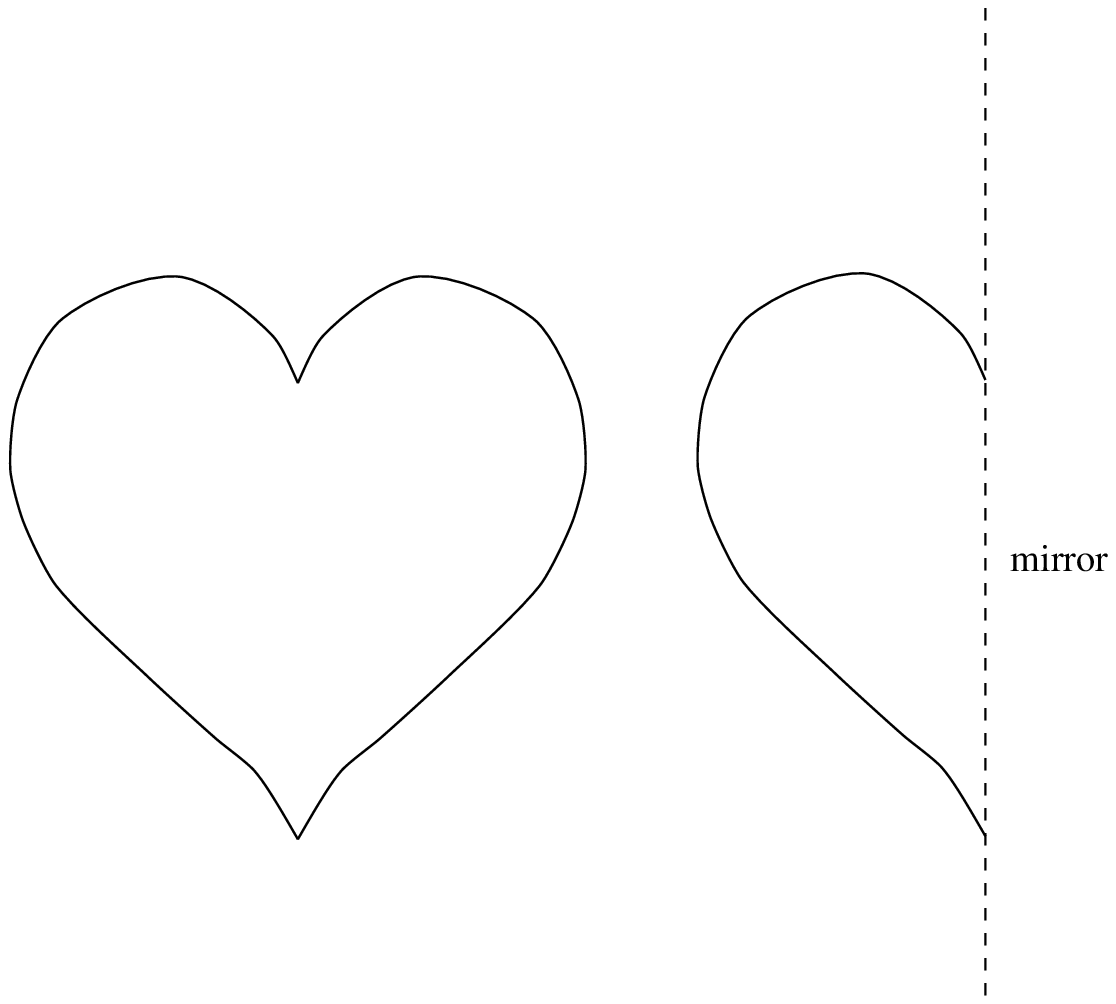}{a heart and its quotient orbifold}{A heart is obtained by
folding a sheet of paper in half, and cutting out half a heart.  The
half-heart is the orbifold for the pattern.  A heart can also be recreated
from a half-heart by holding it up to a mirror.}

A second instance is the paper doll pattern.  Here, there are two
different fold lines.  You make paper dolls by folding a strip of paper
zig-zag, and then cutting out half a person.  The half-person is enough
to reconstruct the whole pattern.  The quotient orbifold is a half-person,
with two mirror lines.

\fig{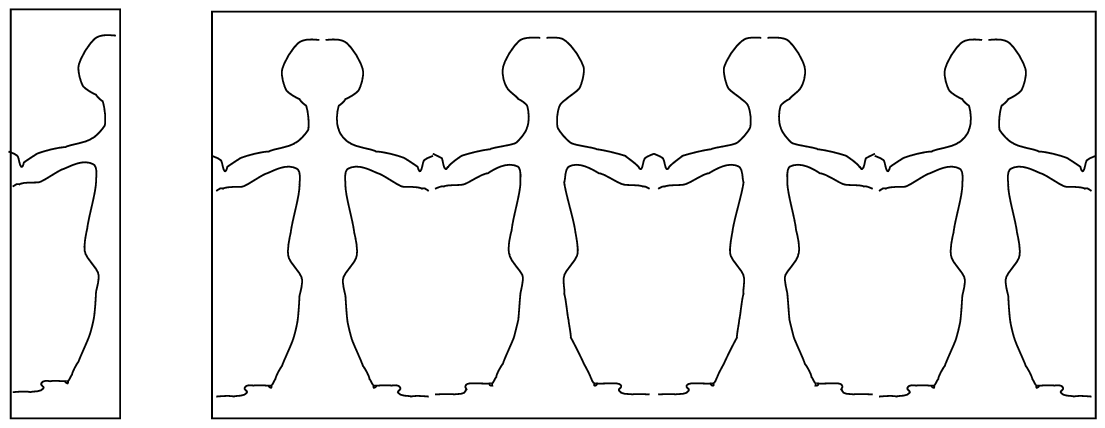}{paper dolls}{A string of paper dolls}

A wave pattern is the next example.  This pattern
repeats horizontally, with no reflections or rotations.  The wave pattern
can be rolled up into a cylinder.  It can be constructed by rolling up
a strip of paper around a cylinder, and cutting a single wave,
through several layers, with
a sharp knife.  When it is unrolled, the bottom part will be like the waves.

\fig{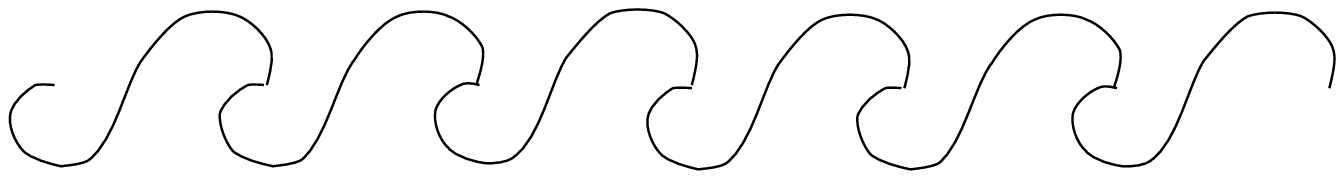}{waves}{This wave pattern repeats horizontally, with no reflections
or rotations.  The quotient orbifold is a cylinder.}

When a pattern repeats both horizontally and vertically, but without
reflections or rotations, the quotient orbifold is a torus.  You can
think of it by first rolling up the pattern in one direction, matching
up equivalent points, to get a long cylinder.  The cylinder has a pattern
which still repeats vertically.  Now coil the cylinder
in the other direction to match up equivalent points on the cylinder.
This gives a torus.

\fig{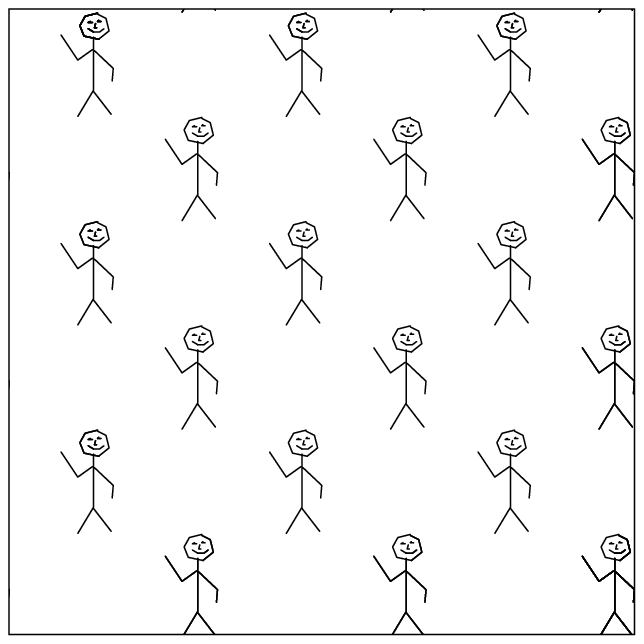}{torus symmetry}{This pattern has quotient orbifold a torus.
It repeats both horizontally and vertically, but without any reflections
or rotations.  It can be rolled up horizontally to form a cylinder,
and then vertically (with a twist) to form a torus.}

\fig{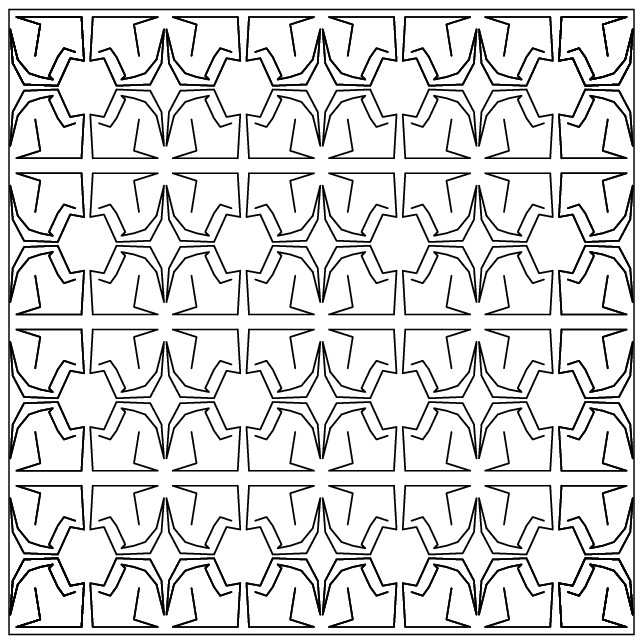}{billiard symmetry}{The quotient orbifold is a rectangle,
with four mirrors around it.}

\fig{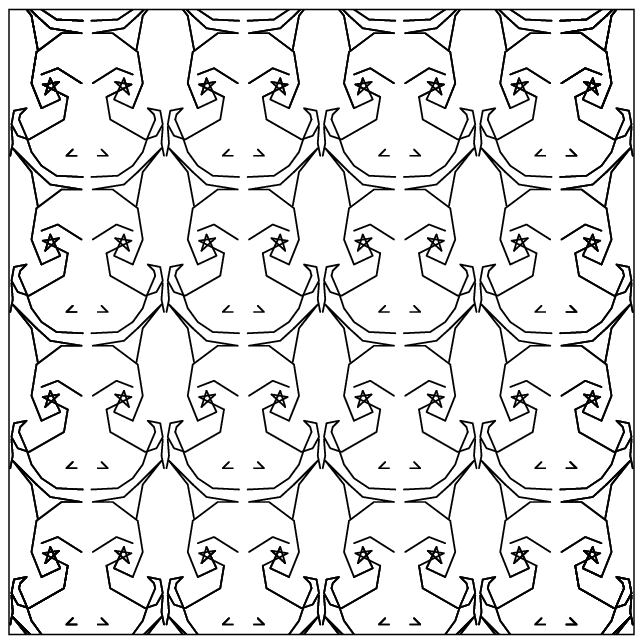}{annulus symmetry}{The quotient orbifold is an annulus,
with two mirrors, one on each boundary.}

\fig{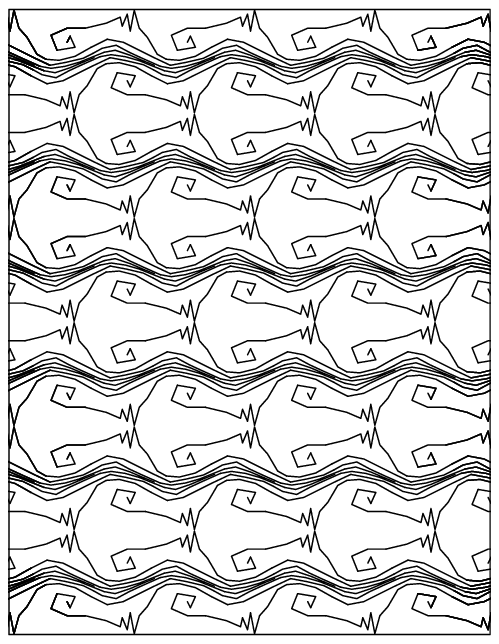}{Moebius symmetry}{The quotient orbifold is a Moebius band,
with a single mirror on its single boundary.}

\fig{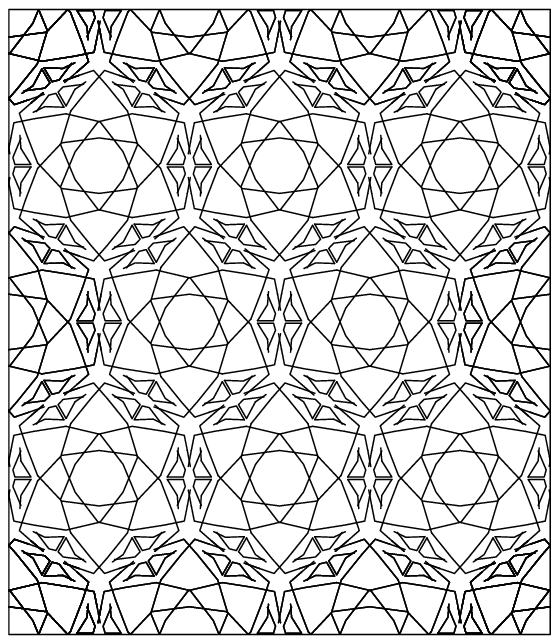}{632 mirrors}{The quotient orbifold is a $60\degree$,
$30\degree$, $90\degree$ triangle, with three mirrors from sides.}

\fig{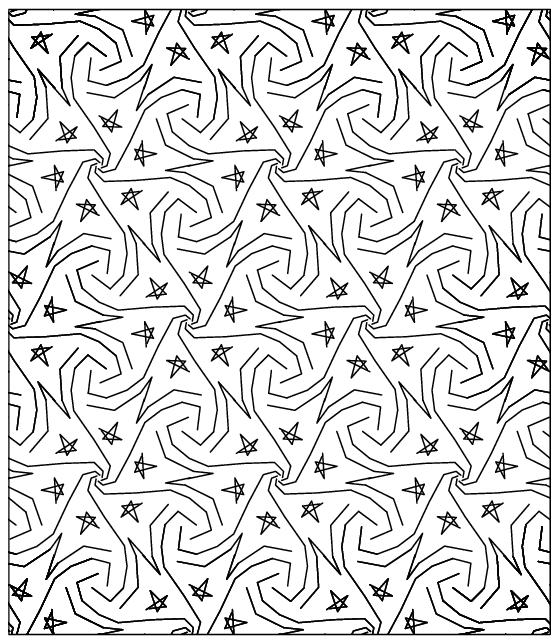}{632 orbifold}{This pattern has rotational symmetry about
various points, but no reflections.  The rotations are of order 6,
3 and 2.  The quotient orbifold is a triangular pillow, with
three cone points.}

\subsection{Discussion}
Using the notation we have discussed, try to figure out the
description of the various pieces of fabric we have handed out.
That is, locate the mirror strings, gyration points, cone points, {\it etc.}
Find the orders of the gyration points and the cone points.

\section{Names for features of symmetrical patterns}

We begin by introducing names for certain features that may occur in symmetrical
patterns.  To each such feature of the pattern,
there is a corresponding feature of the quotient orbifold,
which we will discuss later.

\subsection{Mirrors and mirror strings}
A \dt{mirror} is a line about which the pattern has mirror symmetry.
Mirrors are perhaps the easiest features to pick out by eye.

At a \dt{crossing point}, where two or more mirrors cross,
the pattern will necessarily also have rotational symmetry.
An $n$-way crossing point is one where precisely $n$ mirrors meet.
At an $n$-way crossing point, adjacent mirrors meet at an angle of $\pi/n$.
(Beware: at a 2-way crossing point,
where two mirrors meet at right angles,
there will be 4 slices of pie coming together.)

We obtain a \dt{mirror string}  by starting somewhere on a mirror
and walking along the mirror to the next crossing point,
turning as far right as we can so as to walk along another mirror,
walking to the next crossing point on it, and so on.
(See figure \ref{billiard orbifold}.)

\figsize{150pt}{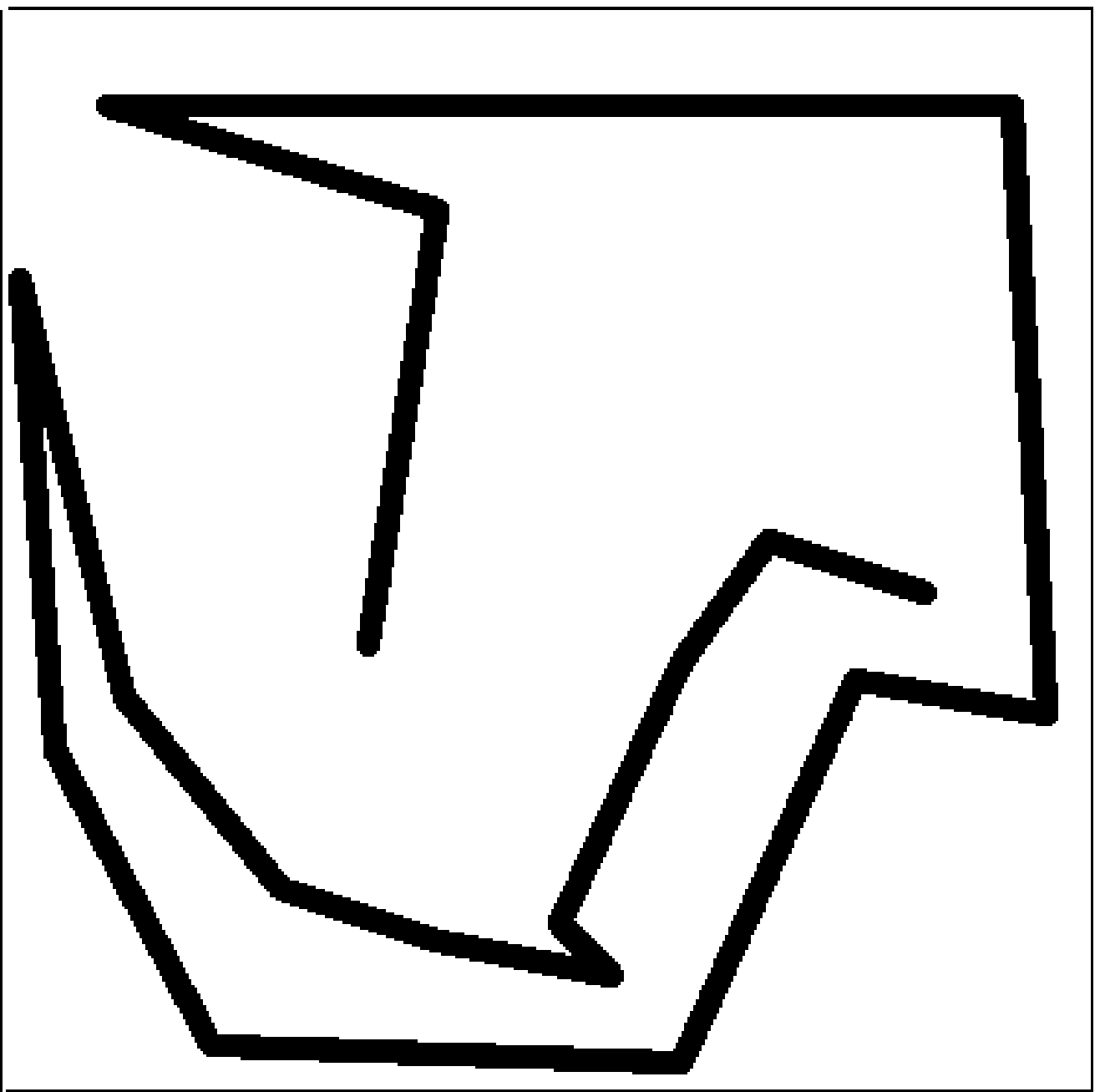}{billiard orbifold}{The quotient billiard orbifold.}

Suppose that you walk along a mirror string
until you first reach a point exactly like the one you started from.
If the crossings you turned at were (say)
a 6-way, then a 3-way, and then a 2-way crossing,
then the mirror string would be of \dt{type} $*632$, etc.
As a special case,
the notation $*$ denotes a mirror that meets no others.

For example, look at a standard brick wall.
There are horizontal mirrors that each bisect a whole row of bricks,
and vertical mirrors that pass through bricks and cement
alternately.
The crossing points, all 2-way, are of two kinds:
one at the center of a brick,
one between bricks.
The mirror strings have four corners,
and you might expect that their type would be $*2222$.
However, the correct type is $*22$.
The reason is that after going only half way round, we come to a point
exactly like our starting point.

\subsection{Mirror boundaries}
In the quotient orbifold, a mirror string of type $*abc$
becomes a boundary wall,
along which there are corners of angles $\pi/a, \pi/b, \pi/c$.
We call this a \dt{mirror boundary} of type $*abc$.
For example, a mirror boundary with no corners at all has type $*$.
The quotient orbifold of a brick wall has a mirror boundary with just
two right-angled corners, type $*22$.

\subsection{Gyration points}
Any point around which a pattern has rotational symmetry is called
a \dt{rotation point}.
Crossing points are rotation points, but there may also be others.
A rotation point that does {\bf \Large NOT} lie on a mirror is called a
\dt{gyration point}.
A gyration point has \dt{order} $n$ if the smallest angle of any rotation
about it is $2 \pi/n$.

For example, on our brick wall there is an order 2 gyration point
in the middle of the rectangle outlined by any mirror string.

\subsection{Cone points}
In the quotient orbifold, a gyration point of order $n$ becomes a cone point
with cone angle $2 \pi / n$.

\section{Names for symmetry groups and orbifolds}

A symmetry group is the collection of all symmetry operations of a pattern.
We give the same names to symmetry groups as to the corresponding
quotient orbifolds.

We regard every orbifold as obtained from a sphere by adding cone-points,
mirror boundaries,
handles, and cross-caps.
The major part of the notation enumerates
the orders of the distinct cone points,
and then the types of all the different mirror boundaries.
An initial black spot $\bullet$ indicates the addition of a handle;
a final circle $\circ$ the addition of a cross cap.

For example, our brick wall gives $2*22$,
corresponding to its gyration point of order 2,
and its mirror string with two 2-way corners.

Here are the types of some of the patterns shown in section 31:

Figure 14: $\bullet$;
Figure 15: $*2222$;
Figure 16: $**$;
Figure 17: $*\circ$.
Figure 18: $*632$.
Figure 19: $632$.

Appart from the spots and circles, these can be read directly
 from the pictures: The important thing to remember is that if
 two things are equivalent by a symmetry, then you only record 
one of them.
A dodecahedron is very like a sphere.
The orbifold corresponding to its symmetry group is a
spherical triangle having angles $\pi/5,\pi/3,\pi/2$;
so its symmetry group is $*532$.

You, the topologically spherical reader,
approximately have symmetry group $*$,
because the quotient orbifold of a sphere by a single reflection
is a hemisphere whose mirror boundary has no corners.

\section{Stereographic Projection}

We let $G$ be a  sphere in Euclidean three space. We want to obtain
a {\it picture} of the sphere on a flat piece of paper or a plane.
Whenever one projects a higher dimensional object onto a lower dimensional object, some type of distortion must occur. There are a number of different ways to project and each projection  preserves some things and distorts others.
 Later  we will explain why we choose {\bf stereographic
 projection}, but first we describe it.

\subsection{Description}
We shall map the  sphere $G$ onto the plane containing its equator. 
 Connect a typical point $P$ on the surface of the sphere to the north pole $N$ by a straight line in three space. This line will intersect the equatorial plane at some point $P'$. We call $P'$ {\bf the projection} of $P$.

 Using this recipe
every point of the sphere except the North pole projects to some point on the equatorial plane. Since we  want to include the North pole in our picture,  we add an extra point $\infty$,  called 
 {\bf the point at infinity},  to the equatorial  plane 
 and we view $\infty$
 as the image of $N$ under stereographic projection.

\subsection{Discussion}
\begin{itemize}
\item Take $G$ to be the unit sphere, $\{(x,y,z) | x^2 + y^2 + z^2 = 1 \}$  so that  $xy$ plane is  the equatorial plane.  The typical point $P$ on the sphere
has coordinates $(X,Y,Z)$. The typical point $P'$ in the equatorial plane, whose coordinates are $(x,y,0)$, will be called $(x,y)$.
 \begin{enumerate}
\item Show that the South pole is mapped into the origin under stereographic projection.
\item Show that under stereographic projection
 the equator is mapped onto the  unit circle, that is the circle $x^2 + y^2 = 1$. 
\item Show that under stereographic projection 
 the lower hemisphere is mapped into the interior of this circle, that is the disk $D = \{ (x,y)| x^2 + y^2 < 1 \}$. 

\item Show that under stereographic projection
 the upper  hemisphere is mapped into the exterior of this circle, that is into $  \{ (x,y)| x^2 + y^2 > 1 \}$. 

 For this to be true where do we have to think of $\infty$ as lying:  interior to $D$ or exterior to it?

\item  
 What projects on to the $x$-axis? 

 What projects onto the $x-\hbox{axis } \cup \infty$?       
Call the set of points that project onto $x-\hbox{axis } \cup \infty$
 the prime meridian. 

\item The prime meridian divides the sphere into two hemispheres, the front hemisphere and the back hemisphere. What is the image of the back hemisphere under
stereographic projection? The front hemisphere?

\item Under stereographic projection what is the image of a great circle passing through the north pole? Of any circle (not necessarily a great circle) passing through the north pole?

\item Under stereographic projection, what projects onto the $y$-axis?
onto any vertical line, not necessarily the $y$ axis?

\end{enumerate}

\end{itemize}
\subsection{What's good about stereographic projection?}
 Stereographic projection preserves circles and angles.  
  That is, the image of a circle on the sphere is a circle in the plane and the angle between two {\it lines} on the sphere is the same as the angle between their images in the plane. A projection that preserves angles is called a {\bf conformal} projection.

We will outline two proofs of the fact that stereographic projection preserves
circles, one algebraic and one geometric. They appear below. 

 Before you do either proof, you may want to clarify in 
your own mind what a {\bf circle} on the surface of a sphere is.
 A circle lying on the sphere is the intersection of a plane in three space with the sphere. 
This can be described algebraically.
For example,  the  sphere of radius 1 with center at the origin is given by
\begin{equation}
G = \{(X,Y,Z)  | X^2 + Y^2 + Z^2 = 1 \}. \label{eq:sphere}
\end{equation}
An arbitrary plane in three-space is given by 

\begin{equation}
AX + BY + CZ + D = 0 \label{eq:plane}
\end{equation}
 for some arbitrary choice of the constants $A$,$B$, 
$C$, and $D$. 
Thus a circle on the unit sphere is any set of points whose coordinates 
simultaneously satisfy equations ~\ref{eq:sphere} and ~\ref{eq:plane}.

\subsubsection{The algebraic proof}
The fact that the  points $P$, $P'$ and $N$ all lie on one line can be expressed by the fact that
\begin{equation}
(X,Y,Z-1) = t(x,y,-1) \label{eq:ratio}
\end{equation}
for some non-zero real number $t$. (Here $P=(X,Y,Z), N=(0,0,1), \hbox { and }
 P'=(x,y,0)$.)

The idea of the proof is that one can 
use   equations ~\ref{eq:sphere} and ~\ref{eq:ratio} to write $X$ as a function of $t$ and $x$, $Y$ as a function of $t$ and $y$, and $Z$ as a function of $t$ and
to simplify equation ~\ref{eq:plane} to an equation in $x$ and $y$.
 Since the equation in $x$ and $y$ so obtained is clearly the equation of a circle in the $xy$ plane, the projection of the intersection of ~\ref{eq:sphere} and ~\ref{eq:plane} is a circle.

To be more precise:

 Equation ~\ref{eq:ratio} says that $X = tx, Y = ty, \hbox{ and } 1-Z = t$.
Set  $Q = \frac{1+Z}{1-Z}$ and verify that
$$  Z = \frac{Q-1}{Q+1},  1+Q = \frac{2}{t}, \hbox { and } Q = x^2 + y^2.$$

 If $P$ lies on the plane, $$AX + BY +CZ + D = 0.$$
Thus
 $$ Atx + Bty  +C\frac{Q-1}{Q+1} + D = 0.$$
Or
$$ \frac{2Ax}{Q+1} + \frac{2By}{Q+1} +C\frac{Q-1}{Q+1} + D = 0.$$
Whence,
$$ 2Ax + 2By + C(Q-1) + D(Q+1) = 0.$$
Or
$$ (C+D)Q + 2Ax + 2By + D-C = 0.$$
Recalling that $Q = x^2 + y^2$, we see
\begin{equation}
  (C+D)(x^2 + y^2) + 2Ax + 2By + D-C = 0 \label{eq:circle}
\end{equation}
 
Since the coefficients of the $x^2$ and the $y^2$ terms are the same, this is 
the equation of a circle in the plane.  

\subsubsection{The geometric proofs}

The geometric proofs sketched below use the following principle:

It doesn't really make much difference if instead of projecting onto
the equatorial plane, we project onto another horizontal plane (not
through N), for example the plane that touches the sphere at the South 
pole, S.  Just what difference does this make?

\begin{itemize}
\item{\bf Angles:}
    To see that stereographic projection preserves angles at $P$, we
project onto the horizontal plane $H$ through $P$.  Then by symmetry
the tangent planes $tN$ and $tP$ at $N$ and $P$ make the same angle $\phi$
with $NP$, as also does $H$, by properties of parallelism (see figure \# 1 at the end of this handout).

  So $tP$ and $H$ are images of each other in the (``mirror'') plane $M$
through $P$ and perpendicular to $NP$.

  For a point $Q$ on the sphere near $P$, the line $NQ$ is nearly parallel to $NP$, so that for points near $P$, stereographic projection is approximately 
the reflection in $M$. 

\item{\bf Circles:} To see
 that stereographic projection takes circles to circles, first
note that any circle $C$ is where some cone touches the sphere, say the cone
of tangent lines to the sphere from a point $V$. 

Now project onto the horizontal plane $H$ through $V$.

In figure \# 2  which {\bf NEED NOT} be a vertical plane,
the four angles $\phi$ are equal, for the same reasons as before, so that
$VP' = VP$. The image of $C$ is therefore the horizontal circle of the
same radius centered at $V$.  

\item {\bf Inversion:}    Another proof uses the fact that
 stereographic projection may be
regarded as a particular case of inversion in three dimensions. You
might like to prove that inversion preserves angles and circularity
in two dimensions.  The {\it inverse} of a point $P$ in the
circle of radius $R$ centered at $O$ is the unique point $P'$ on the
ray $OP$ for which  $OP.OP' = R^2$.

\end{itemize}

\section{The orbifold shop}

The Orbifold Shop has gone into the business of installing orbifold parts.
They offer a special promotional deal: a free coupon
for \$2.00 worth of parts, installation included,
to anyone acquiring a new orbifold.

There are only a few kinds of features
for two-dimensional orbifolds, but they can be used in interesting
combinations.
\begin{itemize}
\item
Handle:  $\$2.00$.
\item
Mirror: $\$1.00$.
\item
Cross-cap:  $\$1.00$.
\item
Order $n$ cone point:  $\$1.00 \times (n-1)/n$.  
\item
Order $n$ corner reflector: $.50 \times (n-1)/n$.
Prerequisite: at least one mirror.  Must specify in mirror and position
in mirror to be installed.
\end{itemize}

With the $\$2.00$ coupon, for example, you could order an orbifold with
four order 2 cone points, costing $\$.50$ each.  Or, you could order
an order 3 cone point costing $\$.66\ldots$, a mirror costing $\$1.00$,
and an order 3 corner reflector costing $\$.33\ldots$.

{\bf Theorem.}
If you exactly spend your coupon at the Orbifold Shop, you will have
a quotient orbifold coming from a symmetrically repeating
pattern in the Euclidean plane with a bounded fundamental domain.
There are exactly $17$ different ways to do this, and corresponding
to the $17$ different symmetrically repeating patterns with bounded
fundamental domain in the Euclidean plane.

\fig{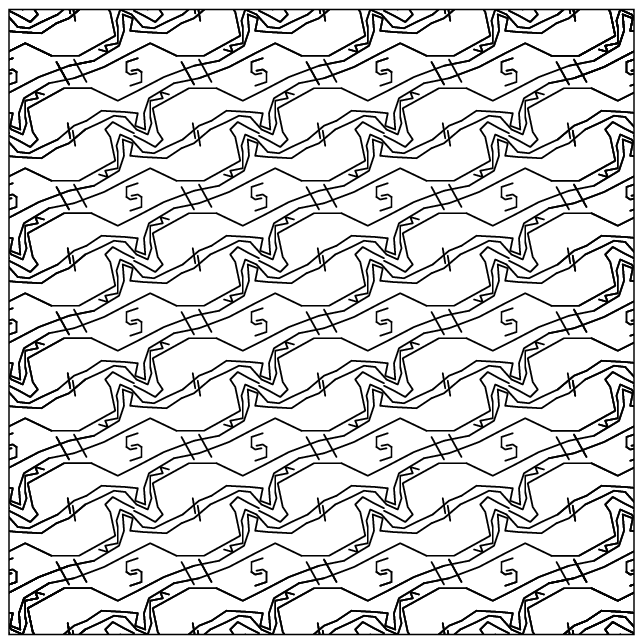}{The 2222 orbifold}{This is the pattern obtained when you
buy four order 2 cone points for $\$.50$ each.}
\fig{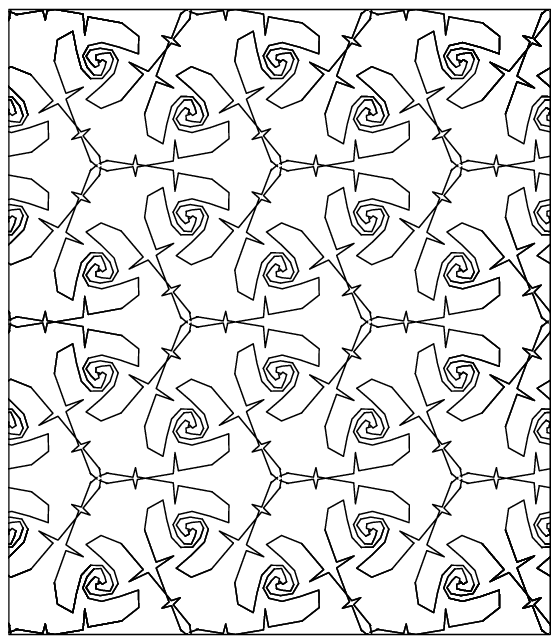}{The 3 star 3 orbifold}{This is the pattern obtained by
buying an order 3 cone point, a mirror, and an order 3 corner reflector.}

{\bf Question.}
What combinations of parts can you find that cost exactly $\$2.00$?

\clearpage
\section{The Euler characteristic of an orbifold}

Suppose we have a symmetric pattern in the plane.  We can make a symmetric
map by subdividing the quotient  orbifold into polygons, and then `unrolling
it' or `unfolding it' to get a map in the plane.

If we look at a large area $A$ in the plane, made up from $N$ copies of a 
fundamental domain, then each face in the map on the quotient orbifold
contributes $N$ faces to the region.  An edge which is not on a mirror
also contributes approximately $N$ copies --- approximately, because when
it is on the boundary of $A$, we don't quite know how to match
it with a fundametnal region.

In general, if an edge or point has order $k$ symmetry which
which preserves it, it contributes approximately $N/k$ copies of
itself to $A$, since each time it occurs, as long as it is not on
the boundary of $A$, it is counted in $k$ copies of the fundamental domain.

Thus,
\begin{itemize}
\item
If an edge is on a mirror, it contributes only approximately
$N/2$ copies.
\item
If a vertex is not on a mirror and not on a cone point, it contributes
approximately $N$ vertices to $A$.
\item
If a vertex is on a cone point of order $m$ it contributes approximately $N/m$
vertices.
\item
If a vertex is on a mirror but not on a corner reflector, it contributes
approximately $N/2$.
\item
If a vertex is on an order $m$ corner reflector, it contributes
approximately $N/2m$
\end{itemize}

{\bf Question.}
Can you justify the use of `approximately' in the list above?  
Take the area $A_R$ to be the union of
all vertices, edges, and faces that intersect a disk of radius $R$
in the plane, along with all edges of any face that intersects and
all vertices of any edge that intersects. Can you show that the
ratio of the true number to the estimated number is arbitrarily close
to 1, for $R$ high enough?

{\bf Definition.}
The \dt{orbifold Euler characteristic} is $V - E + F$, where each 
vertex and edge is given weight $1/k$, where $k$ is the order of symmetry
which preserves it.

It is important to keep in mind the distinction between the topological
Euler characteristic and the orbifold Euler characteristic.  For instance,
consider the billiard table orbifold, which is just a rectangle.
In the orbifold Euler characteristic,
the four corners each count $1/4$, the four edges count $-1/2$, and the 
face counts 1, for a total of 0.  In contrast, the
topological Euler characteristic
is $4 - 4 + 1 = 1$.

{\bf Theorem.}
The quotient orbifold of for any symmetry pattern in the Euclidean plane
which has a bounded fundamental region has orbifold Euler number 0.

{\bf Sketch of proof:}  take a large area in the plane that is topologically
a disk.  Its Euler characteristic is 1.  This is approximately
equal to $N$ times the orbifold Euler characteristic, for some large $N$,
so the orbifold Euler characteristic must be 0.

\bigskip
How do the people at The Orbifold Shop figure its prices?
The cost is based on the orbifold Euler characteristic: it costs
$\$1.00$ to lower the orbifold Euler characteristic by 1.
When they install a fancy new part, they calculate the difference between
the new part and the part that was traded in.

 For instance, to install
a cone point, they remove an ordinary point. An ordinary point counts 1,
while an order $k$ cone point counts $1/k$, so the difference is $(k-1)/k$.

To install a handle, they arrange a map on the original orbifold so that it
has a square face.  They remove the square, and identify opposite edges
of it.  This identifies all four vertices to a single vertex.  The net
effect is to remove 1 face, remove 2 edges (since 4 are reduced to 2),
and to remove 3 vertices.  The effect on the orbifold Euler characteristic
is to subtract $1 - 2 + 3 = 2$, so the cost is $\$2.00$.

{\bf Question.}
Check the validity of the costs charged by The Orbifold Shop
for the other parts of an orbifold.

To complete the connection between orbifold Euler characteristic and
symmetry patterns, we would have to verify that each of the possible
configurations of parts with orbifold Euler characteristic 0 actually
does come from a symmetry pattern in the plane.  This can be done in
a straightforward way by explicit constructions.  It is illuminating
to see a few representative examples, but it is not very illuminating
to see the entire exercise unless you go through it yourself.

\section{Positive and negative Euler characteristic}

A symmetry pattern on the sphere always gives rise to a quotient orbifold
with positive Euler characteristic.  In fact, if the order of symmetry is
$k$, then the Euler characteristic of the quotient orbifold is $2/k$,
since the Euler characteristic of the sphere is 2.

However, the converse is not true.  Not every collection of parts costing
less than $\$2.00$ can be put together to make a viable pattern for
symmetry on the sphere.  Fortunately, the experts at The Orbifold Shop
know the four bad configurations which are too skimpy to be viable:

\begin{itemize}
\item
A single cone point, with no other part, is bad.
\item
Two cone points, with no other parts, is a bad configuration unless they
have the same order.
\item
A mirror with a single corner reflector, and no other parts, is bad.
\item
A mirror with only two corner reflectors, and no other parts, is bad
unless they have the same order.
\end{itemize}

All other configurations are good.  If they form an orbifold with positive
orbifold Euler characteristic, they come from a pattern of symmetry on
the sphere.

\bigskip
The situation for negative orbifold Euler characteristic is straightforward,
but we will not prove it:

{\bf Theorem.}
Every orbifold with negative orbifold Euler characteristic comes from
a pattern of symmetry in the hyperbolic plane with bounded fundamental
domain.  Every pattern of symmetry in the hyperbolic plane with compact
fundamental domain gives rise to a quotient orbifold with negative
orbifold Euler characteristic.

Since you can spend as much as you want, there are an infinite number of these.

\section{Hyperbolic Geometry}
When we tried to make a closed polyhedron by snapping  together
 seven equilateral triangles so that there were seven at every 
vertex, we were unable to do so. 
Those who persisted and  continued to snap together seven
triangles at each vertex,  actually constructed an approximate
 model of the hyperbolic plane. It is this bumpy sheet with angular excesses all over the place that you
might  think of when you try to visualize the hyperbolic plane.
Since we know that angular excess corresponds to negative curvature, we
see that 
the hyperbolic plane is a negatively curved space.

Hyperbolic geometry is also known as {\bf Non-Euclidean geometry}.
 The latter name reflects the fact that it was originally  discovered by
mathematicians seeking a geometry which failed to satisfy Euclid's parallel
postulate. (The parallel postulate states that through  any point
not on a given line there is precisely one line 
 that does not intersect the given line.) While we will outline the details
of Non-Euclidean geometry and prove that it fails to satisfy the parallel postulate, our main emphasis will be on the {\it feel} of the hyperbolic 
plane and hyperbolic 3-space. 

\subsection{Defining the hyperbolic plane}
There are a number of different models for the hyperbolic plane.
They are, of course, all equivalent. As with any instance when 
there are several ways to describe something, each description has
both advantages and disadvantages. We will describe two
models, the {\bf upper half-plane} model, which we denote by {\bf U}
 and the {\bf unit disc}
model, which we initially denote by {\bf D}.

 It will  generally be clear from the context which
model we are using. Although we first present the upper
 half-plane model and prove most of the  fundamental facts there, we will
 generally
after that  use the unit disc.

\subsubsection{The Upper Half-Plane}
Remember that  the image of
the back hemisphere under stereographic projection is the set of all 
points in the $xy$ plane whose $y$ is positive.  This is the upper half-plane. The prime
meridian projects onto the line $y=0$ to which we have added 
 the point at infinity. We think of the image of the prime meridian as the boundary of the upper half-plane. The line $y=0$ could be referred to as the $x$ axis. 
We will also refer to it as the {\bf real axis}, $R$. 
\fig{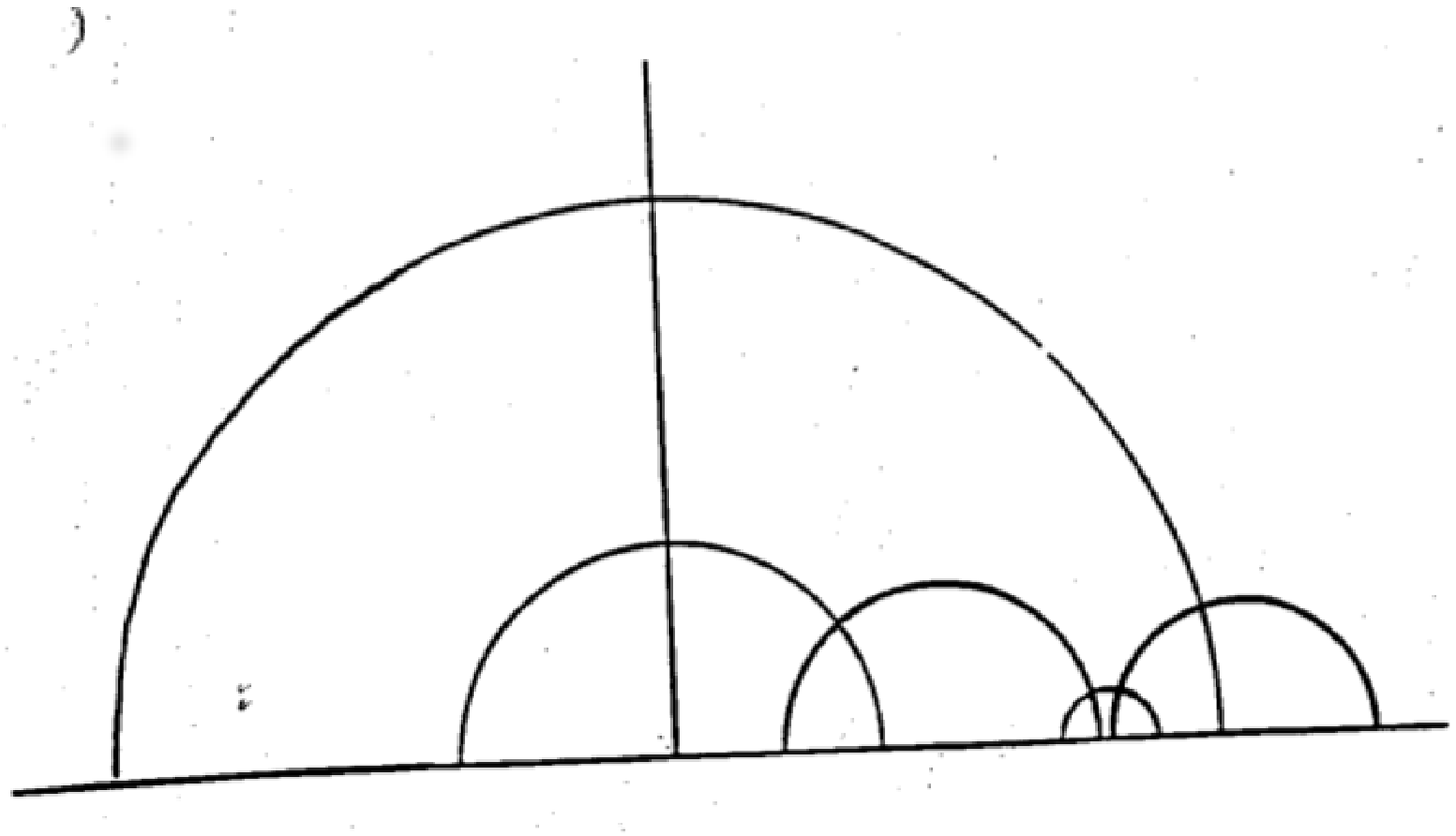}{figure:HU}{Some h-lines in the upper half-plane}
To define a geometry in $U$ we need to define what is meant by a
straight line through two points. 
Given two points $z_1$ and $z_2$ in $U$, one can construct many circles
passing through both of them since {\it three} points determine a circle.
However, there is a unique circle 
passing through $z_1$ and $z_2$ that is perpendicular to $R$. We call
 this circle the {\it straight line} passing through $z_1$ and $z_2$.
 When we want to 
emphasize that we are talking about the hyperbolic line through
two points rather than the Euclidean line, we refer to it as an {\bf h-line}.

You will notice that if $z_1$ and $z_2$ lie on a vertical line, then
there is no Euclidean circle through both that is perpendicular to the boundary. 
However,   recall that   a circle
 on the sphere that passes through the north pole projects onto a what at first glance looks like a line, but is upon reflection can be viewed as a circle (it passes through infinity where it closes up.).      We {\it shall} also view this line as a circle. 

Thus an h-line is either a circle perpendicular to the real axis
 or a vertical line (see figure ~\ref{figure:HU}).
 (The latter is also automatically perpendicular to the real axis.)

(As a homework exercise you can remind yourself
that any circle that intersects the real axis at right angles
has its center on the real axis.)

We note that two Euclidean circles are either disjoint, intersect in a point,
or intersect in two points. Two circles whose centers are on the
real axis that intersect in two points have one point of intersection
above the real axis and one below. Thus they have only one point of
intersection in the upper half-plane. Similarly, a circle with center
on the real axis and a vertical line can have at
 most one point of intersection in the upper half-plane.
Thus any two h-lines are either disjoint or intersect in a point.
We have now proved that this system of {\it lines} and {\it points}
satsifies two of the axioms for a geometry.          
\fig{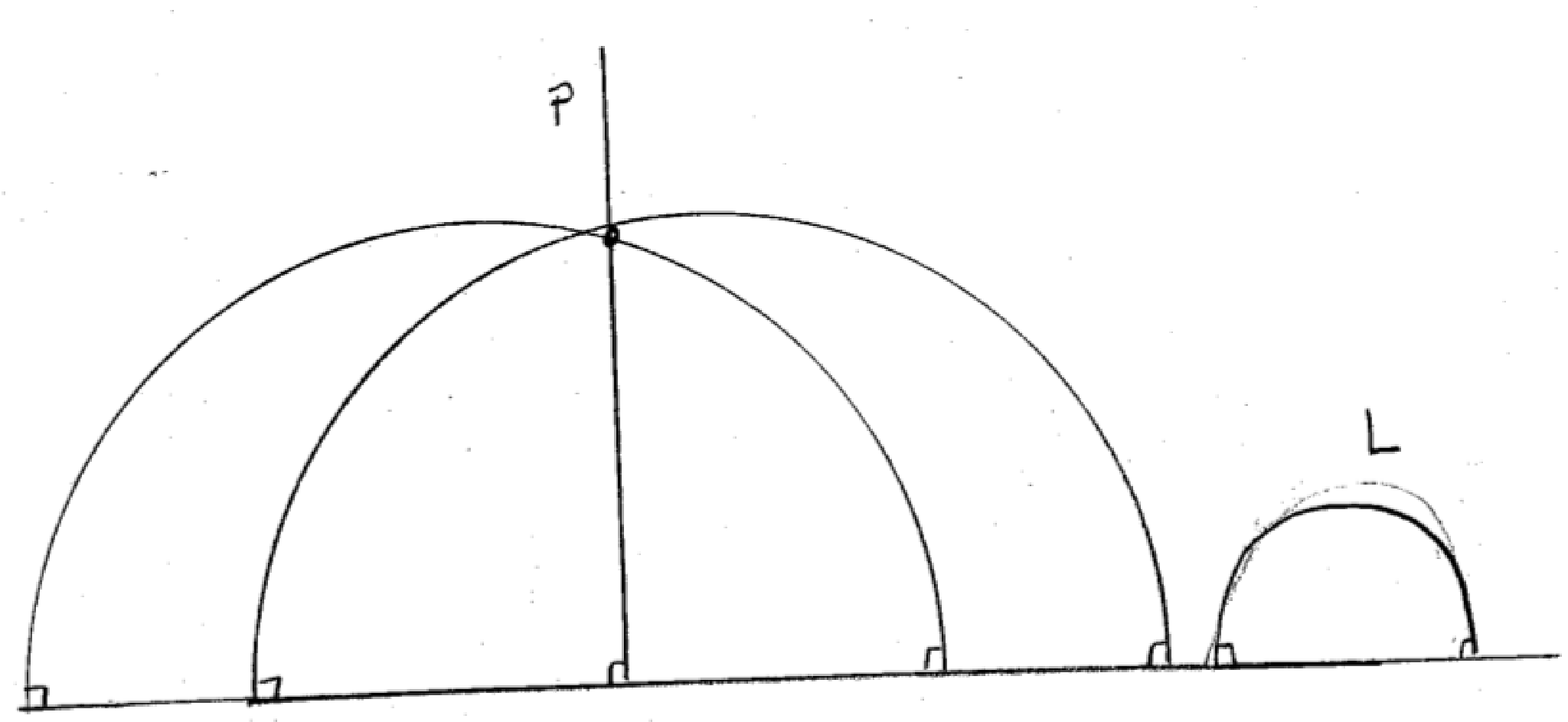}{figure:family}{Several h-lines through $p$ that are disjoint from $L$.}
\subsection{Discussion}
\begin{itemize}
\item What does a hyperbolic mirror look like?
\item  What does a hyperbolic mirror
string look like?
\item What is the maximum number of mirrors in a Euclidean mirror string?
\item What is the maximum number of mirrors in a hyperbolic mirror string?
\end{itemize}
We turn to the parallel axiom. Again let $L$ be any h-line. For the sake of 
simplicity assume that $L$ is not a vertical h-line. Let $p$ be any point
not on $L$. 
We can construct a whole family of Euclidean circles whose centers
are on the real axis which pass through $p$ and which do not intersect
$L$. Figure  ~\ref{figure:family} illustrates several such h-lines. 
For homework, you can work out either one example or a detailed proof. 

\subsection{Distance}
We have emphasized that one of the main distinctions between
geometry and topology is that distance is intrinsic to geometry.
Thus it behooves us to define a distance in the hyperbolic plane. Again, our
emphasis should not be on computing distance, but on having a
feel for hyperbolic distance. The important fact to remember is: 
\begin{itemize}
\item 

{\large Line segments that appear to be of very different lengths
to our Euclidean eyes may be of the same length when we wear 
hyperbolic glasses and vice-versa.}
\end{itemize}
\section{Distance recipe}
Here is a technical definition of how to compute distance.

\begin{indent} 
Begin with any two points. If $L$ is the $h$-line on which they lie, 
let $L'$ be the line on the back hemisphere that projects onto $L$.
Rotate the sphere so that one of the end points of $L'$ moves to the north
pole, $N$. 
 $L'$ rotates into a new line $L''$ passing through $N$. The projection
of $L''$ is now a vertcal line, $K$. The points $a$ and $b$ have been
 lifted to $L'$ rotated to $L''$ and then projected onto $K$. They are now
called $a'$ and $b'$. We can take the  ratio of the heights of $a'$ and $b'$.
This is almost a distance. However,  distance should  be
symmetric. The ratio of the heights depends upon which point we name first.
Therefore, we take the absolute value of the natural log of the ratio of the heights to be
the distance between $a$ and $b$. 
\end{indent}

\subsection{Examples of distances}

Consider the two pairs of points 
\begin{itemize}
\item $A = (0,4)$ and $B = (0,8)$.
\item $C = (0,8)$ and $D = (0,16)$.
\end{itemize}
To our Euclidean eyes it appears to us that $C$ and $D$ are twice as far
apart as $A$ and $B$. When we put on our hyperbolic glasses,
we realize that the distance between $A$ and $B$ is exactly the same
 as the distance between $C$ and $D$.

\subsection{The Unit Disc Model}
 Let $D$ be the unit disc in the plane. $D = \{ (x,y) | x^2 + y^2 < 1 \}$.
We saw earlier that $D$ is the image of the lower half sphere under 
stereographic projection. This is another model for the hyperbolic plane.
We will easily locate the h-lines once we see how this is related to
the upper half-plane.  
\fig{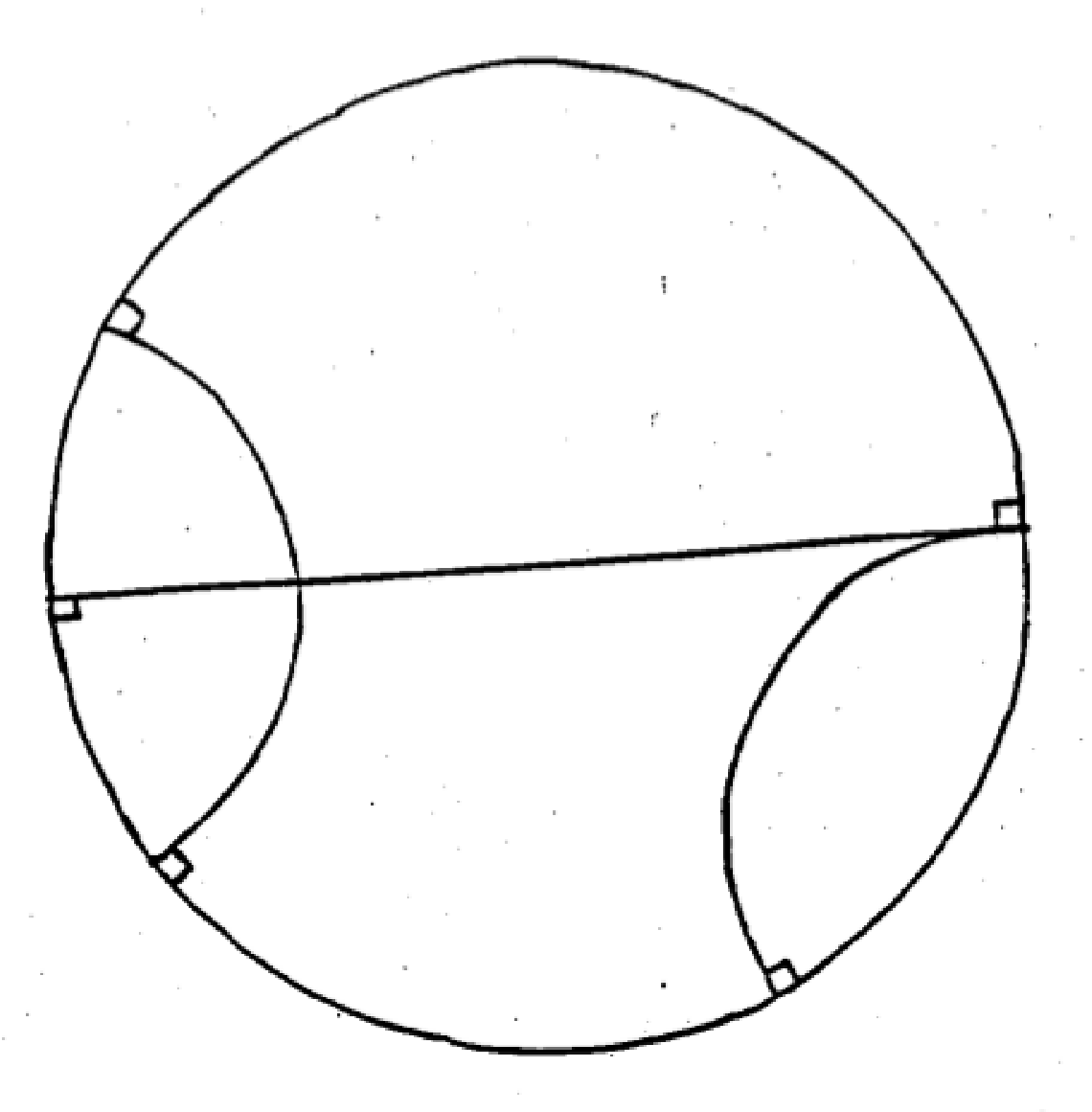}{figure:D}{Some h-lines in $D$.}
\subsection{Passing from one model to another}
Take the sphere. Rotate it so that the back hemisphere goes
 into the bottom hemisphere. Project the bottom hemisphere onto the unit
disc. This procedure identifies the upper half plane (the image of the
back hemisphere) with the unit disc (the image of the bottom hemisphere).
An h-line in the upper half plane corresponds to a circle on the back
hemisphere which is perpendicular to the prime meridian. Such a circle rotates
into a circle on the bottom hemisphere that is perpendicular to the
equator, and then projects to a circle in the plane that intersects the 
boundary of the unit disc at right angles.
 When we project onto the unit
disc, we no longer have to worry about h-lines through infinity. 
Things look much more symmetric. However, we still have one weird
type of h-line: a  Euclidean straight line passing through the center of the
disc. (See figure ~\ref{figure:D}.)
\fig{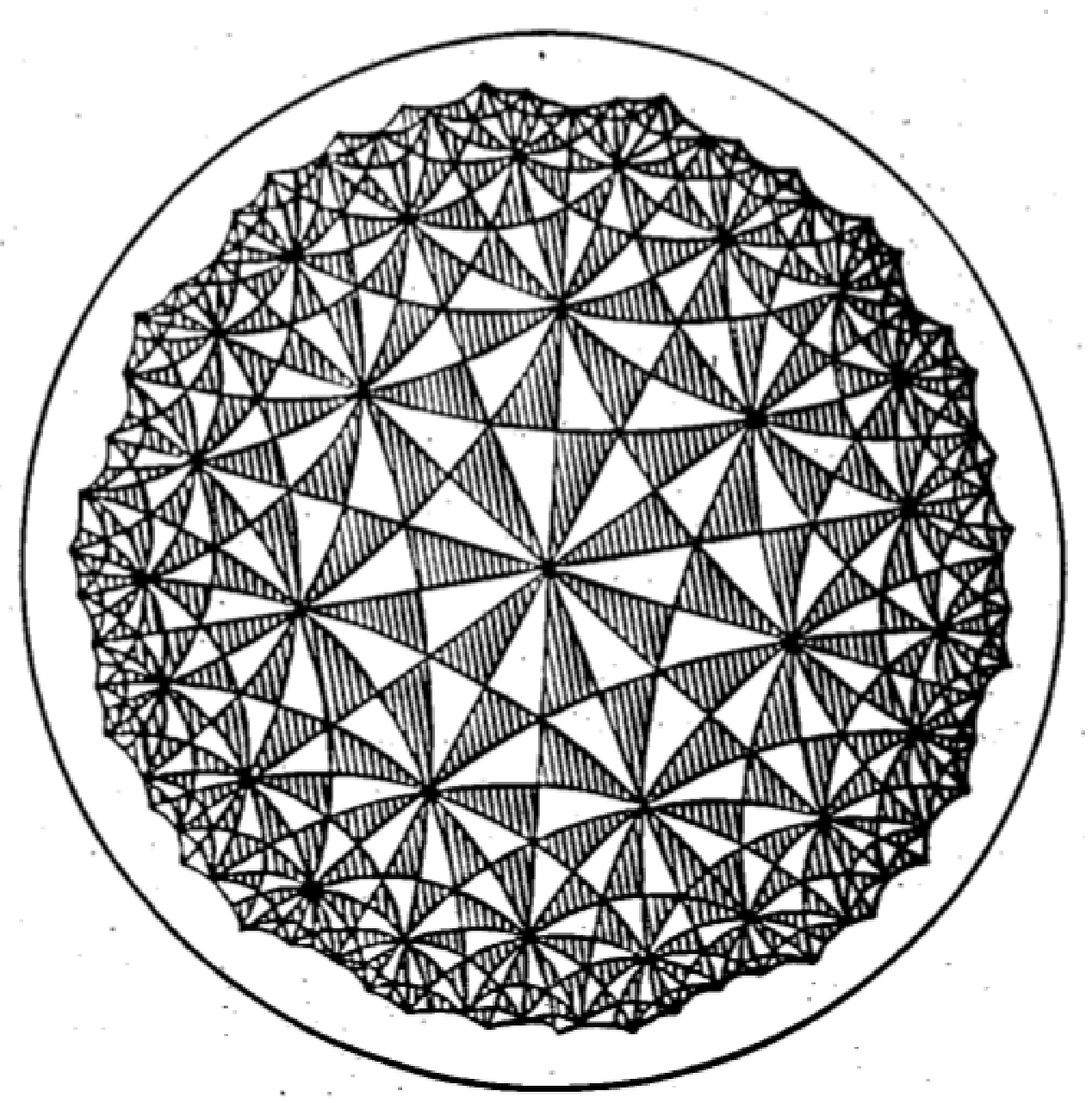}{figure:237}{Some hyperbolic cloth:  A tiling of the hyperbolic plane by triangles with angles $\pi/2,\pi/3,\hbox { and } \pi/7$.}

Once we have a hyperbolic geometry, many new things are possible.
\begin{itemize}
\item  We can classify patterns on hyperbolic cloth. 
We can look for hyperbolic mirrors, hyperbolic gyration points, etc.
and analyze hyperbolic cloth just as we analyzed Euclidean cloth.
\item We can form a $ *237$ orbifold. 
\end{itemize}

Enclosed is a picture of the tiling of the  hyperbolic plane by triangles
whose angles are $\pi/2$, $\pi/3$ and $\pi/7$. (See figure ~\ref{figure:237}.)
The important thing to realize
about this picture is that {\bf ALL} the trianglular tiles  are congruent. That is,
even though the triangles near the boundary of $D$ appear to be much
smaller than those in the center, their sides all have the same lengths.
To see this you just have to look through  your hyperbolic glasses.
\clearpage

\section{A field guide to the orbifolds}
The number 17 is {\em just right} for the number of types of
symmetry patterns in the Euclidean plane:  neither too large nor too small.
It's large enough to make learning to recognize them a challenge,
but not so large that this is an impossible task.
It is by no means necessary to learn to distinguish the 17
types of patterns
quickly,
but if you learn to do it,
it will give you a real feeling of accomplishment,
and it is a great way to amaze and overawe your friends,
at least if they're a bunch of nerds and geeks.

In this section,
we will give some hints about how to learn to classify the patterns.
However,
we want to emphasize that this is a tricky business,
and the only way to learn it is by hard work.
As usual, when you analyze a pattern,
you should look first for the mirror strings.
The information in this section is meant as a way that you
can learn to become more familiar with the 17 
types of patterns,
in a way that will help you to distinguish between them more quickly,
and perhaps in some cases to be able to classify some of the more
complicated patterns without seeing clearly and precisely what the
quotient is.
This kind of superficial knowledge is no substitute for a real
visceral understanding of what the quotient orbifold is,
and in every case you should go on and try to understand why
the pattern is what you say it is while your friends are busy
admiring your cleverness.

This information presented in this section has been gleaned from
a cryptic manuscript discovered among the
personal papers of John Conway after his death.
For each of the 17 types of patterns, the manuscript shows a small
piece of the pattern,
the notation for the quotient orbifold,
and Conway's idiosyncratic pidgin-Greek name for the corresponding
pattern.
These names are far from standard,
and while they are unlikely ever to enter common use,
we have found from our own experience that they are not wholly
useless as a method for recognizing the patterns.

We will begin by discussing Conway's names for the orbifolds.
A reproduction of Conway's manuscript appears at the end of the
section.
You should refer to the reproduction as you try to understand the
basis for the names.

\subsection{Conway's names}
Each of Conway's 17 names consists
of two parts,
a {\em prefix} and a {\em descriptor}.

\subsubsection{The prefix}
The prefix
tells the number of directions from which you can
view the pattern
without noticing any difference.
The possibilities for the prefix are:
{\em hexa-; tetra-; tri-; di-; mono-}.

For example, if you are looking at a standard brick wall,
it will look essentially the same whether you stand on your feet or
on your head.
This will be true even if the courses of bricks in the wall 
do not run parallel to the ground,
as they invariably do.
Thus you can recognize right away that the brick-wall pattern
is {\em di-}something-or-other
In fact, it is {\em dirhombic}.

Another way to think about this is that if you could manage to
turn the brick wall upside down, you wouldn't notice the difference.
Again, this would be true even if you kept your head tilted to one side.
More to the point, try looking at a dirhombic pattern drawn on a sheet
of paper.
Place the paper at an arbitrary angle,
note what the pattern looks like in the large,
and rotate the pattern around until it looks in the large like it did
to begin with.
When this happens, you will have turned the paper through half a rev.
No matter how the pattern is tilted originally, there is always one
and only one other direction from which it appear the same
in the large.

This `in the large' business means that you are not supposed to notice
if, after twisting the paper around,
the pattern appears to have been shifted by a translation.
You don't have to go grubbing around looking
for some pesky little point about which to rotate the pattern.
Just take the wide, relaxed view.

\subsubsection{The descriptor}
The descriptor represents an attempt on Conway's part
to unite patterns that
seem more like each other than they do like the other patterns.
The possibilities for the descriptor are:
{\em scopic; tropic; gyro; glide; rhombic}.

The {\em scopic} patterns are those that emerge from kaleidoscopes:
$*632=$ {\em hexascopic};
$*442=$ {\em tetrascopic};
$*333=$ {\em triscopic};
$*2222=$ {\em discopic};
$**=$ {\em monoscopic};

Their $*$-less counterparts
are the {\em tropic} patterns (from the Greek for `turn'):
$632=$ {\em hexatropic};
$442=$ {\em tetratropic};
$333=$ {\em tritropic};
$2222=$ {\em ditropic};
$\bullet=$ {\em monotropic}.

With the scopic patterns, it's all done with mirrors,
while with the tropic patterns,
it's all done with gyration points.
The two exceptions are:
$**=$ {\em monoscopic};
$\bullet=$ {\em monotropic}.
There is evidence that Conway did not consider these
to be exceptions, on the grounds that
`with the scopics it's all done with mirrors and translations,
while with the tropics,
it's all done with turnings and translations'.

The {\em gyro} patterns contain both mirrors and gyration points:
$4*2=$ {\em tetragyro};
$3*3=$ {\em trigyro};
$22*=$ {\em digyro}.

Since both tropic and gyro patterns involve gyration points,
there is a real possibility of confusing the names.
Strangely, it is the tropic patterns that are the more closely
connected to gyration points.
In practice, it seems to be easy enough to draw this distinction
correctly,
probably because the tropics correspond closely to the scopics,
and `tropic' rhymes with `scopic'.
Conway's view appears to have been that a gyration point,
which is a point of rotational symmetry that does {\bf \Large NOT}
lie on a mirror,
becomes ever so much more of a gyration point when there are mirrors
around that it might have been tempted to lie on,
and that therefore patterns that contain both gyration points and
mirrors are more gyro than patterns with gyration points but no
mirrors.

The {\em glide} patterns involve glide-reflections:
$22\circ=$ {\em diglide};
$\circ\circ=$ {\em monoglide}.

The glide patterns are the hardest to recognize.
The quotient orbifold of the diglide pattern is a projective plane with
two cone points;
the quotient of the monoglide patterns is a Klein bottle.
When you run up against one of these patterns,
you just have to sweat it out.
One trick is that when you meet something that has glide-reflections
but not much else,
then you decide that it must be either a diglide or a monoglide,
and you can distinguish between them by deciding whether it's a {\em di-}
or a {\em mono-} pattern,
which is a distinction that is relatively easy to make.
Another clue to help distinguish these two cases
is that a diglide pattern has glides in two different
directions, while a monoglide has glides in only one direction.
Yet another clue is that in a monoglide you can often spot two disjoint
M\"{o}bius strips within the quotient orbifold,
corresponding to the fact that the quotient orbifold for a monoglide
pattern is a Klein bottle, which can be pieced together from two M\"{o}bius
strips.
These two disjoint M\"{o}bius strips arise from the action
of glide-reflections along parallel but inequivalent axes.

The {\em rhombic} patterns often give a feeling of rhombosity:
$2*22=$ {\em dirhombic};
$*\circ=$ {\em monorhombic}.

An ordinary brick wall is dirhombic;
it can be made monorhombic by breaking the 
gyrational symmetry.
The quotient of a monorhombic pattern is a M\"{o}bius strip.
Like the two glide quotients, it is non-orientable,
but it is much easier to identify
because of the presence of the mirrors.

\subsection{How to learn to recognize the patterns}
As you will see, Conway's manuscript
shows only a small portion of each of the patterns.
A very worthwhile way of becoming acquainted with the patterns is to
draw larger portions of the patterns,
and then go through and analyze each one,
to see why it has the stated notation and name.
You may wish to make flashcards to practice with.
When you use these flashcards,
you should make sure that you can not only remember the correct
notation and name, but also that you can analyze the pattern quickly,
locating the distinguishing features.
This is important because the patterns you will see in the real world
won't be precisely these ones.

Another hint is to keep your eyes open for symmetrical patterns in the
world around you.
When you see a pattern,
copy it onto a flashcard,
even if you cannot analyze it immediately.
When you have determined the correct analysis,
write it on the back and add it to your deck.

\subsection{The manuscript}
What follows is an exact reproduction of Conway's manuscript.
In addition to the 17 types of repeating patterns,
Conway's manuscript also gives tables of the 7 types of
frieze patterns, and of the 14 types of symmetrical patterns on the sphere.
These parts of the manuscript appear to be mainly gibberish.
We reproduce these tables here in the hope that they may someday come to the
attention of a scholar who will be able to make sense of them.

\newcommand{\bigfigj}[1]{
\begin{figure}
\centerline{\mbox{\includegraphics[width=400pt]{figures/#1.eps}}}
\end{figure}
}
\bigfigj{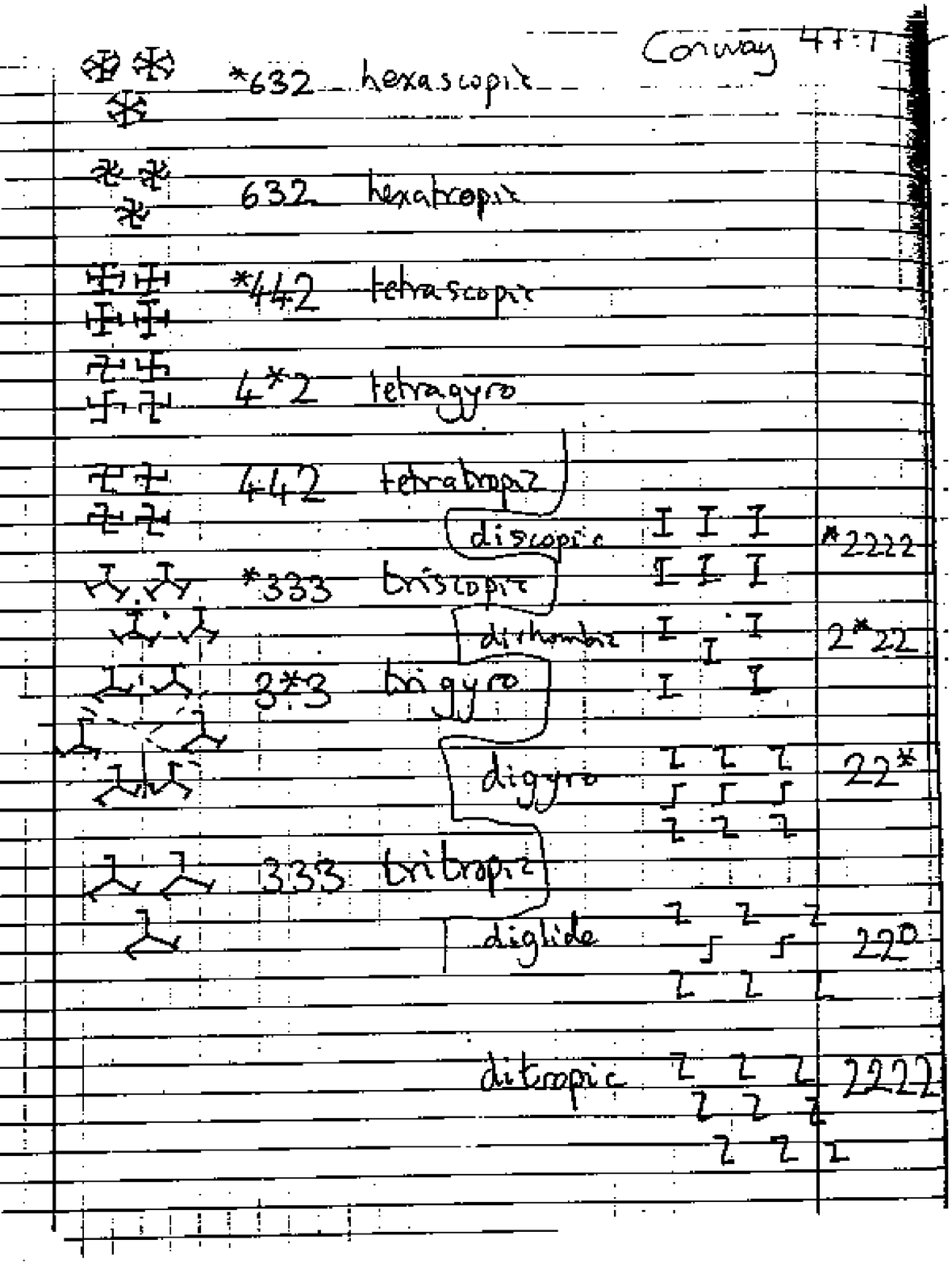}

\bigfigj{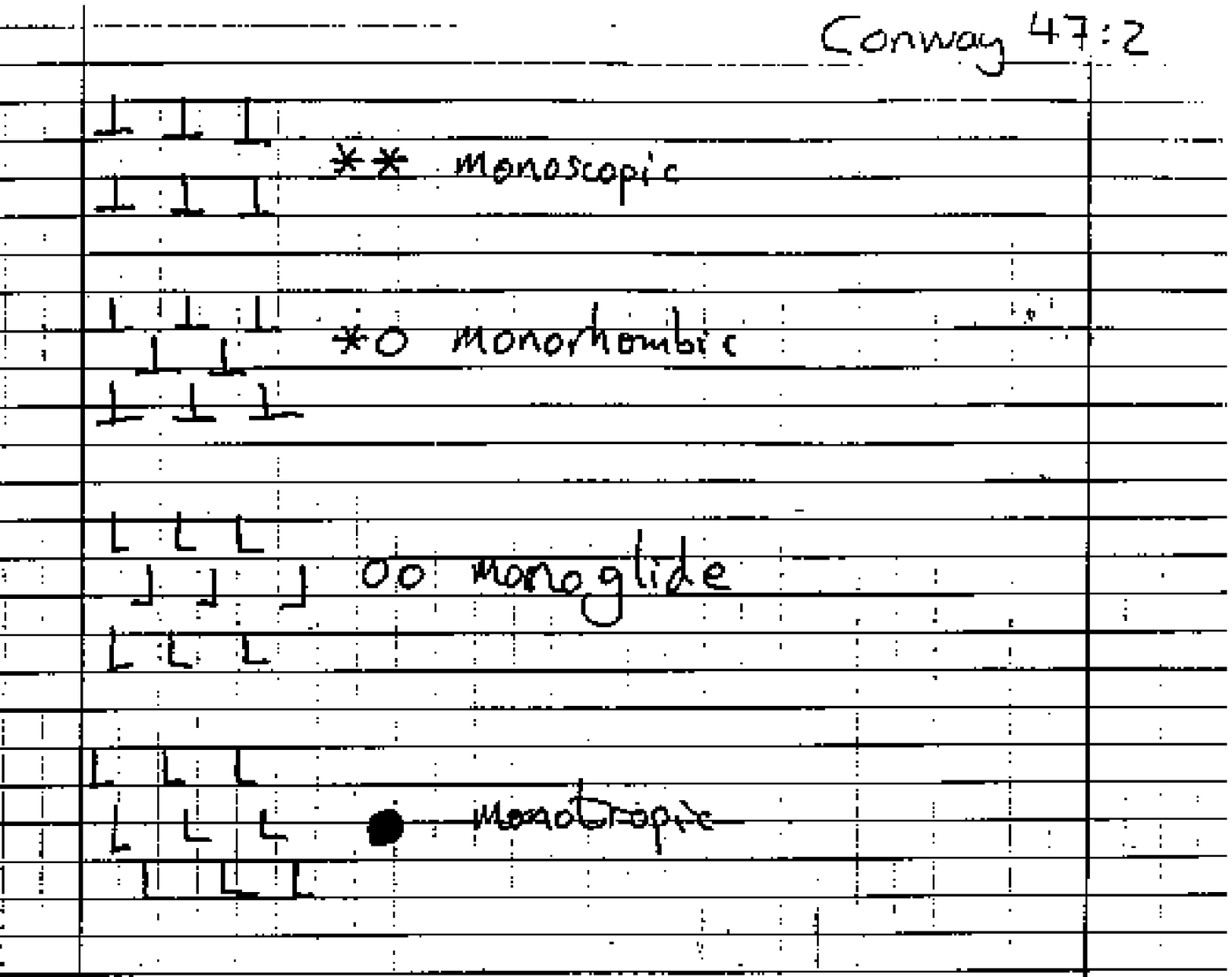}

\bigfigj{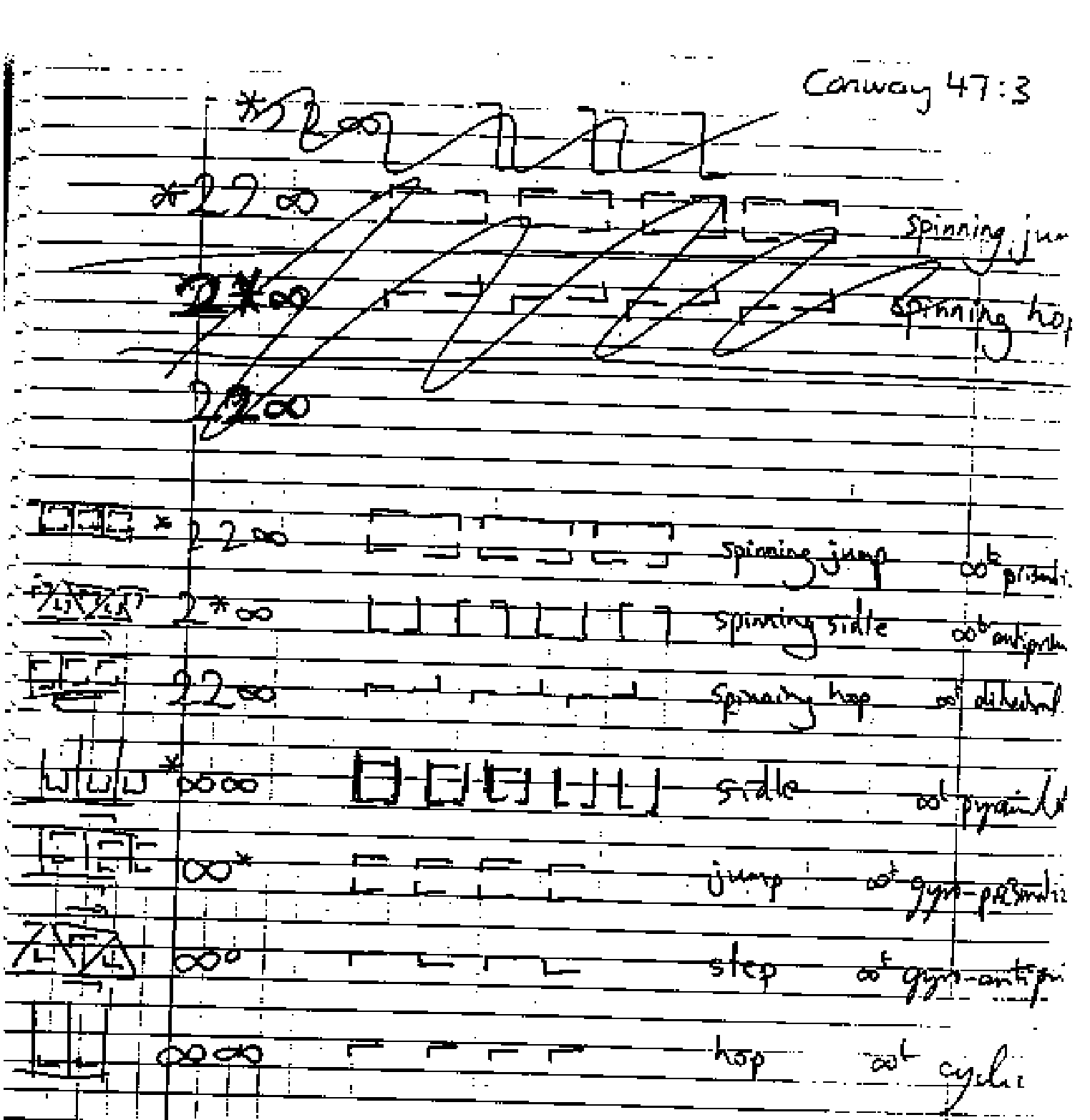}

\bigfigj{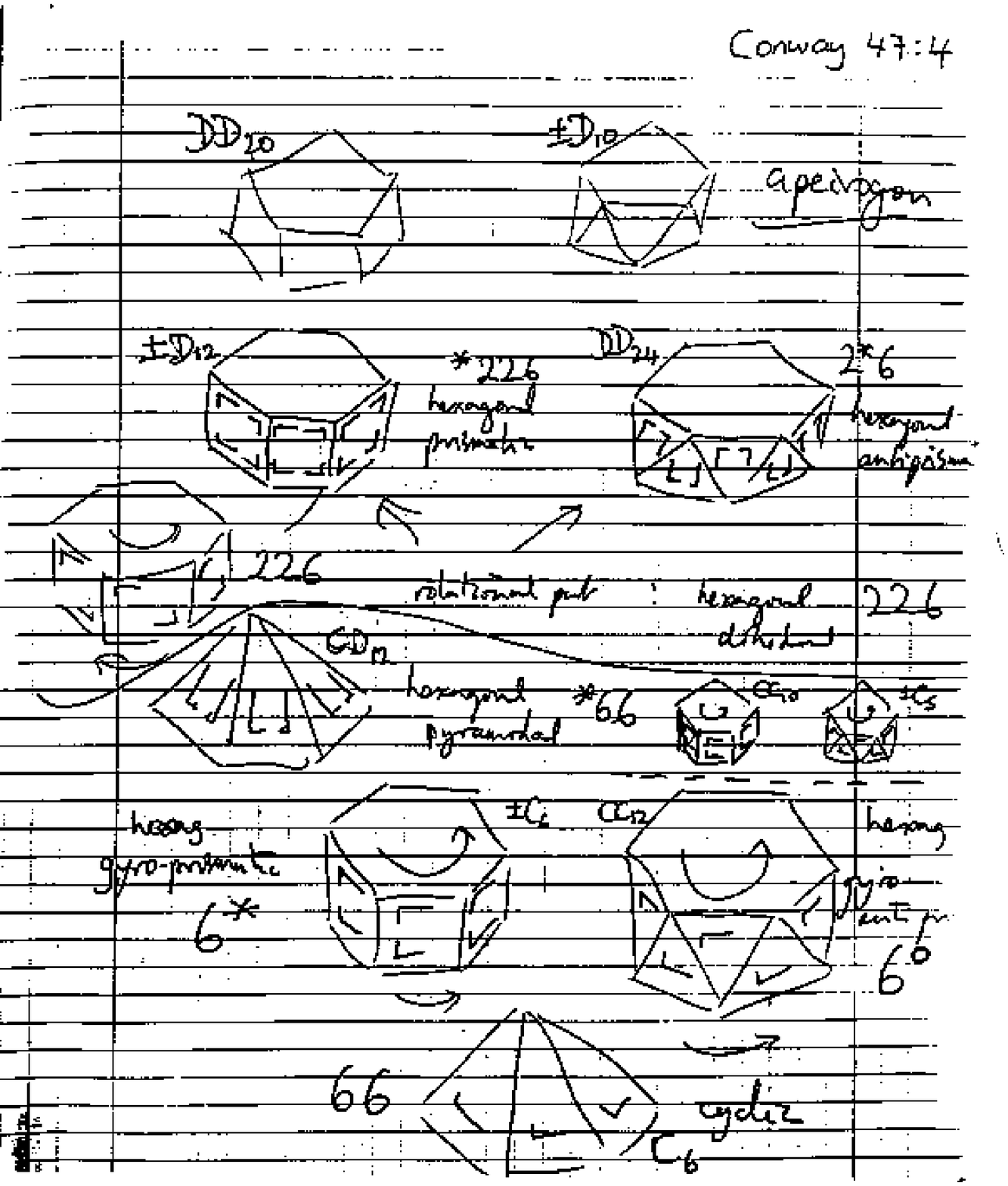}

\bigfigj{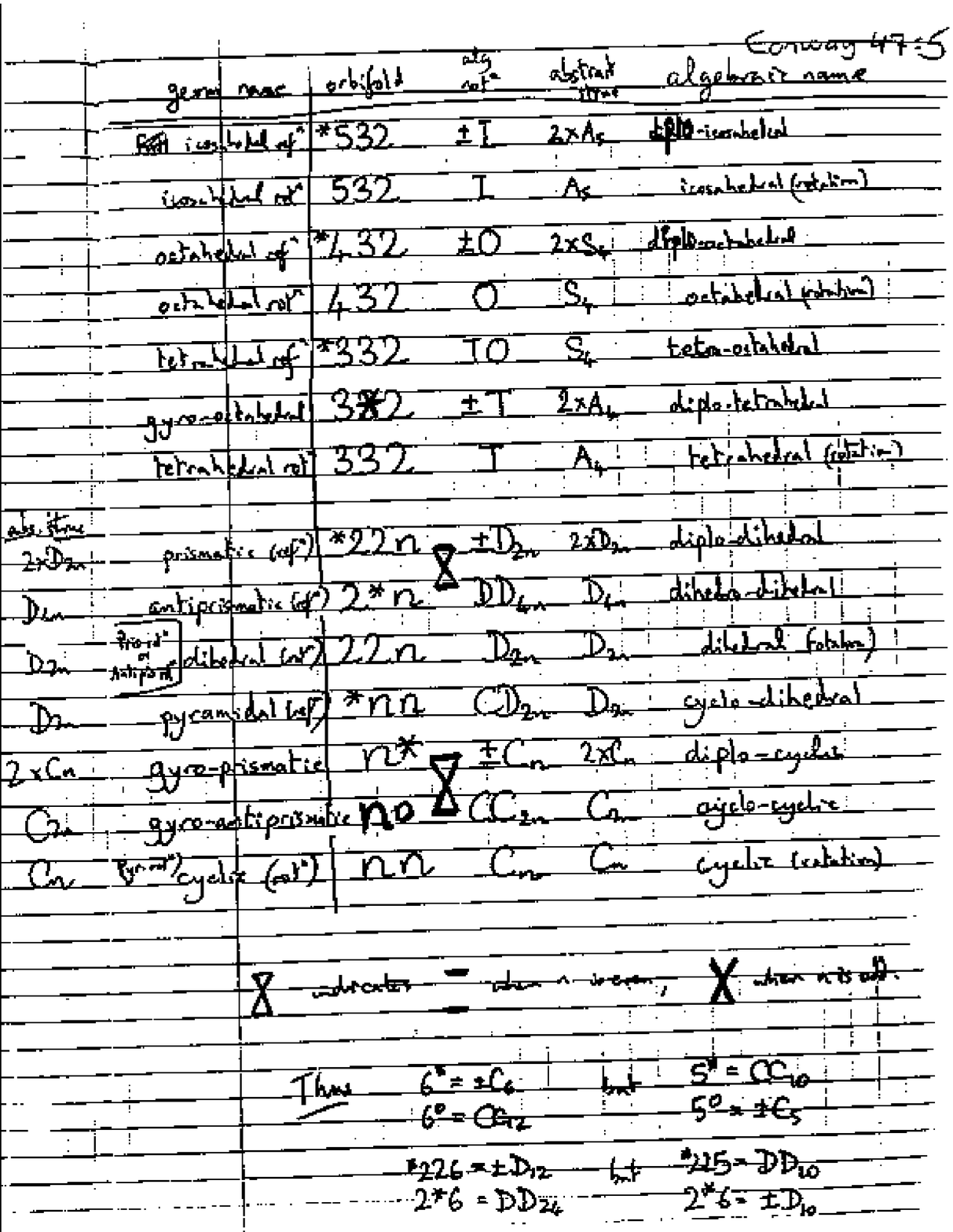}

\end{document}